\documentclass[a4paper,12pt]{amsart}

\usepackage{amssymb}
\usepackage{amsmath}
\usepackage{amsthm}
\usepackage{calc}
\usepackage{latexsym}
\usepackage{mathrsfs}
\usepackage{upgreek}

\allowdisplaybreaks
\numberwithin{equation}{section}

\theoremstyle{plain}
\newtheorem{thm}{Theorem}[section]
\newtheorem{prop}[thm]{Proposition}

\newtheorem{cor}[thm]{Corollary}
\newtheorem*{tom}{\rubrik}
\newcommand{\rubrik}{}

\theoremstyle{definition}
\newtheorem{defn}[thm]{Definition}
\newtheorem{rem}[thm]{Remark}

\newcommand{\R}{\mathbb{R}}

\newcommand{\Z}{\mathbb{Z}}

\newcommand{\calF}{\mathcal{F}}

\newcommand{\calL}{\mathcal{L}}
\newcommand{\calM}{\mathcal{M}}

\newcommand{\calS}{\mathscr{S}}
\newcommand{\calE}{\mathscr{E}}
\newcommand{\FT}{\mathscr F}
\newcommand{\MT}{\mathscr M}
\newcommand\Rn{{{\mathbb R}^n}}
\newcommand\Tn{{{\mathbb T}^n}}

\newcommand{\scal}[2]{\langle #1,#2\rangle}
\newcommand{\rr}[1]{\mathbb R^{#1}}

\newcommand{\nm}[2]{\Vert #1\Vert _{#2}}

\newcommand{\sets}[2]{\{ \, #1\, ;\, #2\, \} }

\newcommand{\fy}{\varphi}
\newcommand{\cdo}{\, \cdot \, }
\newcommand{\supp}{\operatorname{supp}}

\newcommand{\rank}{\operatorname{rank}}

\newcommand{\vrum}{\vspace{0.1cm}}

\def\p#1{{\left({#1}\right)}}
\def\b#1{{\left\{{#1}\right\}}}

\def\n#1{{\left\|{#1}\right\|}}
\def\abs#1{{\left|{#1}\right|}}

\begin{document}

\title{Changes of variables in modulation and Wiener
amalgam spaces}

\author[M. Ruzhansky]{Michael Ruzhansky}

\address{Department of Mathematics, Imperial College London, UK}

\email{m.ruzhansky@imperial.ac.uk}

\author[M. Sugimoto]{Mitsuru Sugimoto}

\address{Department of Mathematics, Graduate School of Science,
Osaka University, Toyonaka, Osaka 560-0043, Japan}

\email{sugimoto@math.sci.osaka-u.ac.jp}

\author[J. Toft]{Joachim Toft}

\address{Department of Mathematics and Systems Engineering,
V{\"a}xj{\"o} University, 351 95 V{\"a}xj{\"o}, Sweden}

\email{joachim.toft@vxu.se}

\author[N. Tomita]{Naohito Tomita}

\address{Department of Mathematics, Graduate School of Science,
Osaka University, Toyonaka, Osaka 560-0043, Japan}

\email{tomita@gaia.math.wani.osaka-u.ac.jp}

\date{\today}

\thanks{The first author was supported by the JSPS Invitational
Research Fellowship.}

%\address{}
%\email{}

%\address{}
%\email{}


%%%==================================================================
\begin{abstract}
In this paper various properties of global and local changes of
variables as well as properties of canonical transforms are
investigated on modulation and Wiener amalgam spaces.
We establish several relations among localisations of
modulation and Wiener amalgam spaces and, as a consequence,
we obtain several versions of local and global Beurling--Helson
type theorems. We also establish a number of positive results
such as local boundedness of canonical transforms on modulation
spaces, properties of homogeneous changes of variables, and
local continuity of Fourier integral operators on $\FT L^q$.
Finally, counterparts of these results are discussed for
spaces on the torus as well as for weighted spaces.

\end{abstract}
%%%==================================================================
\maketitle

%%%==================================================================
%%%==================================================================
\section{Introduction}

The main purpose of this paper is to investigate the invariance
properties of modulation spaces and certain types of Wiener amalgam
spaces under changes of variables. We
%%  will   %%%%
establish different positive
and negative results in these spaces as well as in closely
related Fourier Lebesgue spaces. Let us point out that a natural
ingredient of our analysis is to consider also the
canonical transforms which are changes
of variables on the Fourier transform side.
The canonical transforms play
an %%   (instead of "a very")
important role in the
analysis of partial differential equations because they
allow to transform operators into each other by
changes of variables on the Fourier
transform side (e.g. \cite{DH72}).
Regularity properties of canonical transforms
are important for various applications, for example
in recent applications to
global smoothing problems for evolution equations
(e.g. \cite{RS06a, RS06b}).

Since the Fourier image of
a modulation space is a
%% instead of "the modulation space is the corresponding"
Wiener amalgam space it is natural to consider invariance
properties of changes of variables and canonical transforms
on both spaces. Another space of interest is the
space $\FT L^q$, $1\leq q\leq\infty$, which is the image of
the Lebesgue space $L^q(\Rn)$ under the Fourier transform.
In fact, when localised in space, this space coincides with
modulation spaces $M^{p,q}$ and Wiener amalgam spaces
$W^{p,q}$, so the question of continuity in $\FT L^q(\Rn)$
is related to the question of continuity in its image
under the Fourier transform, which is the usual $L^q(\Rn)$.
For example, when investigating a property of the local
boundedness
of canonical transforms in $L^q(\Rn)$, we can reduce the
analysis to an equivalent question of the Fourier-local boundedness
of changes of variables in $\FT L^q(\Rn)$.
We note that these questions are usually quite delicate since
there is a loss of regularity of Fourier integral operators
in $L^q$--spaces (cf. \cite{SSS91}), which is dependent
on the underlying geometry (cf. \cite{Ru00}).

The question of the invariance of function spaces under changes
of variables is of fundamental importance since it allows to
introduce counterparts of these spaces on manifolds via
localisations. Thus, both local and global invariance properties
are of importance. Unfortunately, many spaces of interest
have a so-called Beurling--Helson property which means that
a $C^1$ change of variables which leaves the space invariant
must be affine (for space $\FT L^1$ on the torus this goes back to
Beurling and Helson \cite{BH}). For example, this property
was established in $\FT L^q$ in \cite{LO94, Se76},
and in modulation spaces in \cite{Ok08}. In Theorem \ref{TH:B-H}
we
%%  will establish it also in %%%%
also establish it for Wiener amalgam spaces.
Our analysis is based on the fact that when
localised in space, function spaces
$M^{p,q}$, $W^{p,q}$ and $\FT L^q$ all coincide
(see Theorem \ref{TH:mw}). This will follow from the
fact that when localised in frequency, function spaces
$M^{p,q}$, $W^{p,q}$ and $L^p$ also coincide.
This observation puts the study of the Beurling--Helson
property on Wiener type spaces in a unified setting,
as well as simplifies the proof in the case
of modulation spaces given in \cite{Ok08}.
In Corollary \ref{cor:spaces} we state various equalities
of localisations of these spaces and Theorem
\ref{TH:B-H}
%%%%will give  %%%%
gives the Beurling--Helson properties
for both changes of variables and canonical transforms.

However, it turns out that we can still prove some
positive results. For example, in Theorem \ref{TH:local-lq}
we will show that if the pullback by a change of variables
$\psi:\Rn\to\Rn$  is bounded on $L^q(\Rn)$ then the corresponding
canonical transform $I_\psi$ (which is the pullback by $\psi$ on the
Fourier transform side) is locally continuous on $M^{p,q}$,
$W^{p,q}$ and $\FT L^q$. On the Fourier transform side this
gives a Fourier-local continuity of the change of
variables induced by such $\psi$
(see Theorem \ref{TH:local-lq} for a
precise statement).

On the other hand, phase functions which come from the
theory of Fourier integral operators are positively
homogeneous of order one (\cite{Hormander}). This means
that the analysis of the invariance properties is important
also outside of the $C^1$ category. In Theorem
\ref{thm:hom} we
%% will %%%%
give a result to this end which
shows that different types of properties are possible.
In particular, we establish a Beurling--Helson type
result in this case as well by using the theory of
Fourier integral operators in an essential way.

At the same time, positive results will allow us to
improve the continuity properties of Fourier integral operators
related to canonical transforms in $\FT L^q$--spaces.
In particular, in \cite{CNR08}, it was shown that
Fourier integral operators are locally bounded on $\FT L^q(\Rn)$
provided that the amplitude is in the symbol
class $S_{1,0}^{-n|1/2-1/q|}$. In
Theorem \ref{cor:FIOs} we remove the decay
condition in the case of canonical transforms and show that
the corresponding operators with amplitudes in
$S^0_{0,0}$ (or even in $M^{\infty,1})$ are still locally
bounded in $\FT L^q(\Rn)$.

Finally, in Theorem \ref{homogen} we investigate other
homogeneous changes of variables which may have singularities
on sets of different dimensions. For them, we
show continuity in modulation and Wiener amalgam spaces. In the
proof of this theorem we use Gabor theory of
modulation spaces and certain decompositions of homogeneous
mappings (cf. Chapter 12 in \cite{Grochenig}). This result
extend previously known properties on $\FT L^q$ and on $M^{p,q}$
with $p=q$.

Modulation spaces were introduced by Feichtinger in \cite{Fe2}
and \cite{Fe4} during the period 1980--1983.
The basic theory of such spaces was thereafter established
and extended by Feichtinger and Gr{\"o}chenig
(see e.{\,}g. \cite{Fe4,FG1,FG2,Grochenig}, and references therein).
Roughly speaking, the (classical) modulation space $M^{p,q}$ is
obtained by imposing a mixed $L^{p,q}$ norm on
the short-time Fourier transform of
a tempered distribution.

\par

A major idea behind these spaces is to find useful
Banach spaces, which are defined in a way similar to Besov spaces, in
the sense of replacing the dyadic decomposition on the Fourier
transform side, characteristic to Besov spaces, with a {uniform}
decomposition. From the construction of these spaces, it turns out
that modulation spaces and Besov spaces in some sense are rather
similar (see \cite{BaC,ST,To04,To04B} for sharp embeddings).

It appears that in some respects these spaces have better
properties from the point of view of evolution partial differential
equations. For example, it was shown in \cite{BGOR} that
propagators for the wave and Schr\"odinger equations are
bounded on modulation spaces, compared to the usual
loss of derivatives in Sobolev spaces
(see e.g. \cite{SSS91}).
We point out in Remark
\ref{rem:propagators} that propagators of the form
$e^{it|D|^\alpha}$ are actually locally continuous on
$M^{p,q}(\Rn)$ and $W^{p,q}(\Rn)$ for all $p,q$ and all
$t, \alpha\in\R$
(compared to the case of $0\leq\alpha\leq 2$ on
$M^{p,q}(\Rn)$ analysed in \cite{BGOR} and to the
well-known loss of
derivatives in local $L^p$ spaces, e.g. for $\alpha=1$ for
the wave equation or
for the KdV equation for $\alpha=3$, etc).

Counterparts of these properties as well as of other results of this
paper for spaces on the torus are discussed in the last section. In
particular, we observe the equality
$M^{p,q}(\Tn)=W^{p,q}(\Tn)=\FT \ell^q(\Tn)$ for all
$1\leq p,q\leq\infty$. This
immediately reduce
the analysis of the Beurling--Helson property to the
original paper of Beurling and Helson \cite{BH} as well
as to the extensions in \cite{LO94}. 
In particular, in the case of $q=1$ the above equality can be 
viewed as a characterisation of 
absolutely convergent Fourier series.

Moreover, we
show the boundedness of canonical transforms on these
spaces. Finally, we will remark that propagators of the form
$e^{it|D|^\alpha}$ are actually isometries on
$M^{p,q}(\Tn)$ for all $p,q$ and all $\alpha\in\R$, and
will discuss periodic weighted spaces.

We note that Theorem \ref{TH:mw} emphasizes difficulties with
the definition of modulation and Wiener amalgam spaces on
manifolds. However, the global definition is still possible in the
presence of the group structure. For example, modulation
spaces on locally compact abelian groups were investigated
in \cite{Fe4}. It is also possible to introduce these spaces
on general compact Lie groups with the global interpretation of
pseudo-differential operators as in \cite{RT-book}.
In this case results of Section 5 can be extended to
the setting of general compact Lie groups.

In Section 2 we state our results. Section 3 will
introduce necessary definitions and terminology.
Proofs and further comments of various nature
will be given in Section 4. Section 5
are devoted to giving some remarks on counterparts of our results
for spaces on the torus. In the appendix we will discuss
weighted spaces.

%%%==================================================================
%%%==================================================================
\section{Results}

First of all we remark some fundamental identities, which
show that in the $\calE'$ and $\FT\calE'$ categories, modulation,
Wiener amalgam, and $\FT L^q$ (or $L^q$) spaces coincide
(the relation between $M^{p,q}$ and $\FT L^q$ spaces has been 
known before, see further
for references). In all sections except for Section 5 we
deal with spaces on $\Rn$.

\begin{thm}\label{TH:mw}
Let $1 \le p,q \le \infty$.
Then the following equalities hold:
\begin{equation}\label{EQ:mw-eq}
\begin{aligned}
&M^{p,q}\cap\calE'=W^{p,q}\cap\calE'=\FT L^q\cap\calE',
\\
&M^{p,q}\cap \FT\calE'=W^{p,q}\cap \FT\calE'=L^p\cap \FT\calE',
\end{aligned}
\end{equation}
with equivalence of norms.
Moreover, let $\Omega \subset \Rn$ be compact. Then the estimates
\begin{equation}\label{EQ:mw-norms}
\begin{aligned}
&\|f\|_{M^{p,q}}\leq C|\widetilde{\Omega}|^{\max(0,1/q-1/p)}
\|f\|_{W^{p,q}},
\\
&\|f\|_{W^{p,q}}\leq C|\widetilde{\Omega}|^{\max(0,1/p-1/q)}
\|f\|_{M^{p,q}}
\end{aligned}
\end{equation}
hold %% added
for all $f \in \calS'$ with $\mathrm{supp}\, f \subset \Omega$,
where constant $C>0$ is independent of $\Omega$ and $f$,
$\widetilde{\Omega}=\{x \in \Rn : \mathrm{dist}(x,\Omega)<1\}$,
and $|\widetilde{\Omega}|$ is the Lebesgue measure of
$\widetilde{\Omega}$.
\end{thm}

We note that \eqref{EQ:mw-norms} is equivalent to
\[
C^{-1}|\widetilde{\Omega}|^{\min(0,1/q-1/p)}\|f\|_{W^{p,q}}
\leq \|f\|_{M^{p,q}}
\leq C|\widetilde{\Omega}|^{\max(0,1/q-1/p)}\|f\|_{W^{p,q}}.
\]
We note also that a weighted version of equalities
\eqref{EQ:mw-eq} will be given in Remark \ref{thmagain}.

%\begin{rem}\label{rem:propagators}
%In \cite{To04} is was proved that $e^{i|D|^2}$ is bounded on each modulation
%space, and in
%\cite{BGOR} it was shown that for $0\leq \alpha\leq 2$,
%operators $e^{i|D|^\alpha}$ are bounded on modulation
%spaces $\calM^{p,q}(\Rn)$, $1\leq p,q\leq\infty$
%(for the definition of $\calM^{p,q}(\Rn)$
%see Remark \ref{modified-def}).
%In particular, this covers wave and Schr\"odinger propagators.
%Theorem \ref{TH:mw} immediately implies that Fourier multipliers
%$e^{i|D|^\alpha}$ are locally continuous on
%$M^{p,q}(\Rn)$ and $W^{p,q}(\Rn)$ for all $1\leq p,q\leq\infty$
%and all $\alpha\in\R$. Indeed, if
%$\chi_1, \chi_2\in C_0^\infty(\Rn)$, then we can estimate
%\begin{multline*}
% \n{\chi_1 e^{i|D|^\alpha}\chi_2 f}_{M^{p,q}(\Rn)}\asymp
% \n{\chi_1 e^{i|D|^\alpha}\chi_2 f}_{\FT L^{q}(\Rn)}=
% \n{\chi_1(D) e^{i|\xi|^\alpha}\chi_2(D)\widehat{f}}_{L^q(\Rn)}
% \\
% \leq C\n{\chi_2(D)\widehat{f}}_{L^q(\Rn)} =
% C\n{\chi_2 f}_{\FT L^q(\Rn)} \asymp
% C\n{\chi_2 f}_{M^{p,q}(\Rn)},
%\end{multline*}
%using that $\chi_1(D)$ is bounded on $L^q(\Rn)$ for all
%$1\leq q\leq\infty$ by the Young's inequality.
%This observation increases the expectation that Fourier
%multipliers $e^{i|D|^\alpha}$ should be bounded on
%$M^{p,q}(\Rn)$ also for $\alpha$ outside of the interval
%$[0,2]$.
%\end{rem}

\begin{rem}\label{rem:propagators}
In \cite{To04} it was proved that $e^{i|D|^2}$ is bounded on each
modulation space, and in
\cite{BGOR} it was shown that for $0\leq \alpha\leq 2$,
operators $e^{i|D|^\alpha}$ are bounded on modulation
spaces $\calM^{p,q}(\Rn)$, $1\leq p,q\leq\infty$
(for the definition of $\calM^{p,q}(\Rn)$
see Remark \ref{modified-def}).
In particular, this covers wave and Schr\"odinger propagators.

\par

On the other hand, Theorem \ref{TH:mw} can be used to establish
local continuity properties for a broader class of Fourier
multipliers. More precisely, assume that $m\in L^\infty (\Rn )$. Then
$m(D)$ from $\mathscr S(\Rn )$ to $\mathscr S'(\Rn )$ extends uniquely
to a  locally continuous map on $M^{p,q}(\Rn)$ and on $W^{p,q}(\Rn)$
for all $1\leq p,q\leq\infty$. Indeed, if  $\chi_1, \chi_2\in
C_0^\infty(\Rn)$, then
\begin{multline*}
 \n{\chi_1 m(D)\chi_2 f}_{M^{p,q}(\Rn)}\asymp
 \n{\chi_1 m(D)\chi_2 f}_{\FT L^{q}(\Rn)}=
 \n{\chi_1(D) m(\xi )\chi_2(D)\widehat{f}}_{L^q(\Rn)}
 \\
 \leq C\nm {m}{L^\infty} \n{\chi_2(D)\widehat{f}}_{L^q(\Rn)} =
 C\nm {m}{L^\infty}\n{\chi_2 f}_{\FT L^q(\Rn)} \asymp
 C\nm {m}{L^\infty}\n{\chi_2 f}_{M^{p,q}(\Rn)},
\end{multline*}
using the fact that $\chi_1(D)$ is bounded on $L^q(\Rn)$ for all
$1\leq q\leq\infty$ by Young's inequality.

\par

In particular we may choose $m(\xi )=e^{i|\xi |^\alpha}$, for any
$\alpha \in \mathbb R$, and this observation together
with the corresponding results on the torus (see
Section 5) increase the expectation
that Fourier multipliers $e^{i|D|^\alpha}$ should be bounded on
$M^{p,q}(\Rn)$ also for $\alpha$ outside of the interval
$[0,2]$.
\end{rem}

Since we are going to investigate properties of operators
in localisations of function spaces both in space and in
frequency, it is convenient to introduce the following
notation. Let $X\subset\calS'(\Rn)$ be a normed linear
spaces. Then we introduce the following notation for
functions which are compactly supported either in space
or in frequency
\begin{equation}\label{EQ:notation1}
X_{comp}:= X\cap \calE', \quad X_{\FT comp}:=X\cap \FT\calE',
\end{equation}
as well as localisations of these spaces
\begin{equation}\label{EQ:notation2}
\begin{aligned}
& X_{loc}:= \{u\in\calS': \chi u\in X \textrm{ for all }
\chi\in C_0^\infty(\Rn)\}, \\
&X_{\FT loc}:= \{u\in\calS': \chi(D) u\in X \textrm{ for all }
\chi\in C_0^\infty(\Rn)\}.
\end{aligned}
\end{equation}
All these spaces inherit the metric from $X$ and
from $\FT X$ in a natural way.
We will say that for normed linear
spaces $X, Y\subset\calS'$, a mapping
$T:X\to Y$ is {\em locally bounded} 
(or locally continuous)
if it is continuous
from $X_{comp}$ to $Y_{loc}$, and that it is
{\em Fourier--locally bounded} 
(or Fourier--locally continuous)
if it is continuous
from $X_{\FT comp}$ to $Y_{\FT loc}$.
%We also say that a mapping $T:X\to Y$ is 
%{\em Fourier--continuous} if 
%$\FT\circ T\circ\FT^{-1}$ is
%continuous from $\FT X$ to $\FT Y$.

\par

Since
$$
M^{p,q}_{comp}=W^{p,q}_{comp}=(\FT L^q)_{comp},\qquad
M^{p,q}_{\FT comp}=W^{p,q}_{\FT comp}=L^p_{\FT comp},
$$
by Theorem \ref{TH:mw}, and since trivially
$$
\FT(M^{p,q}_{comp})=W^{q,p}_{\FT comp},\qquad
\FT(W^{p,q}_{comp})=M^{q,p}_{\FT comp}
$$
we obtain the following corollary to Theorem \ref{TH:mw}:

\begin{cor}\label{cor:spaces}
Let $1 \le p,q \le \infty$. Then the equalities
\begin{alignat*}{4}
\FT(M^{p,q}_{comp}) &= W^{q,p}_{\FT comp} & &= L^q_{\FT comp} & &=
M^{q,p}_{\FT comp} & &=\FT(W^{p,q}_{comp})
\intertext{and}
\FT(M^{p,q}_{\FT comp}) &= W^{q,p}_{comp} & &= (\FT L^p)_{comp} & &=
M^{q,p}_{comp} & &=\FT(W^{p,q}_{\FT comp})
\end{alignat*}
hold.
\end{cor}

Now we introduce two important operators.
Given a mapping $\psi$ from $\R^n$ to itself, we
define the change of variables $\psi^*$ by
\[
(\psi^*f)(x)=f(\psi(x))
\]
and the canonical transform $I_\psi$ by
\[
I_\psi f(x)=\FT^{-1}[(\FT f)(\psi(\xi))](x)
\]
for functions $f$ on $\R^n$.
Clearly, we have the equality
\begin{equation}\label{EQ:ch-can}
I_\psi=\FT^{-1}\circ\psi^*\circ\FT.
\end{equation}
The combination of Theorem \ref{TH:mw} and known Beurling--Helson
type theorems
give   
the following Beurling--Helson
local and global type theorems for modulation and Wiener
amalgam spaces:

\begin{thm}\label{TH:B-H}
Let $1 \le p,q \le \infty$, $ 2 \neq q<\infty$,
and let $\psi : \R^n \to \R^n$ be a $C^1$-function.
Assume that one of the following conditions are fulfilled:
\begin{itemize}
\item[(i)] operator $\psi^*$
is bounded on either $M^{p,q}(\R^n)$, $W^{p,q}(\R^n)$ or
 $\FT L^q(\Rn)$;
\item[(ii)] operator $\psi^*$
is locally bounded on either $M^{p,q}(\R^n)$, $W^{p,q}(\R^n)$ or
 $\FT L^q(\Rn)$;
\item[(iii)] operator $I_\psi$
is bounded on either $M^{q,p}(\R^n)$, $W^{q,p}(\R^n)$ or $L^q(\Rn)$;
\item[(iv)] operator $I_\psi$
is Fourier--locally
bounded on either $M^{q,p}(\R^n)$, $W^{q,p}(\R^n)$ or $L^q(\Rn)$;
\end{itemize}
Then $\psi$ is an affine function.
\end{thm}

It was pointed out in \cite{Ok08}
that in the case of $M^{p,q}$ in condition (i)
the statement essentially reduces to the Beurling--Helson
type theorem on $\FT L^q$ which was treated earlier in
\cite{BH, LO94, Se76}. We will give a simplified proof of such reduction
using Theorem \ref{TH:mw}, with a simple proof of equalities
\eqref{EQ:mw-eq}, at least in the case when one does not need
to keep track of constants in \eqref{EQ:mw-norms}.
Theorem \ref{TH:mw} also allows us to treat the Wiener
amalgam spaces (so we formulate Theorem \ref{TH:B-H}
in a unified way).

We note that pairs of assumptions (i)-(ii) and (iii)-(iv)
in Theorem \ref{TH:B-H}
are obviously equivalent in view of \eqref{EQ:ch-can}.
However, we choose to write all of them explicitly because of
the following result that shows that
non-affine transforms can be allowed if we
just consider the
\emph{local} boundedness of $I_\psi$ on modulation spaces,
or if we localise a
change of variables on the Fourier transform side
in Wiener amalgam spaces:

\begin{thm}\label{TH:local-lq}
Let $1 \le p,q \le \infty$,
and let $\psi : \R^n \to \R^n$ be such that
$\psi^*$ is bounded on $L^q(\Rn)$.
Then the following is true:
\begin{itemize}
\item[(i)]
$I_\psi$ is locally continuous on $M^{p,q}(\Rn)$,
$W^{p,q}(\Rn)$ and $\FT L^q(\Rn)$;
\item[(ii)] $\psi^*$ is Fourier-locally continuous on
$M^{p,q}(\Rn)$, $W^{p,q}(\Rn)$ and $L^p(\Rn)$.
\end{itemize}
\end{thm}

Another important class of canonical transforms that arises
in applications to partial differential equations and in
the theory of Fourier integral operators is the
class of functions $\psi$ positively homogeneous of
order one, which means that $\psi(\lambda x)=\lambda\psi(x)$
for all $\lambda>0$ and all $x\in\Rn$. In this case,
this function is no longer $C^1$ everywhere, and we have
mixed results already on the space $\FT L^q(\Rn)$:

\begin{thm}\label{thm:hom}
Let $\psi:\Rn\to\Rn$ be positively
homogeneous of order one and let $q$ be such that
$1\leq q\leq\infty$. Then the following is true:
\begin{itemize}
\item[(i)] assume that the inverse $\psi^{-1}$ exists
on $\Rn\backslash 0$ and
satisfies
$\psi^{-1}\in C^{1}(\Rn\backslash 0)$. Then
$\psi^*$ is Fourier-locally continuous on $L^q(\Rn)$;
\item[(ii)]
assume that $\psi\in C^\infty(\Rn\backslash 0)$. Assume
also that $\psi^*$ is continuous or
Fourier-locally continuous on
$\FT L^q(\Rn)$ and $q\not=2$. Then $\psi$ is linear.
\end{itemize}
\end{thm}
%{\bf Question: what about if $\psi^*$ is locally continuous
%on $\FT L^q(\Rn)$? In fact, we are ultimately interested whether
%$I_\psi$ is bounded on $M^{p,q}$. If this is true, then
%$I_\psi$ is also Fourier-locally bounded on $L^p$
%(i.e. bounded on $L^p_{\FT loc}$), which would imply that
%$\psi^*$ is locally bounded on $\FT L^p(\Rn)$.}

%\bigskip

If $q=2$, then clearly $\psi^*$ is continuous
(and hence also Fourier-locally continuous) on
$\FT L^2(\Rn)=L^2(\Rn)$.
By using relation \eqref{EQ:ch-can} we can easily obtain
a counterpart of this theorem for canonical transforms $I_\psi$.
We note that part (i) is a straightforward consequence of
Theorem \ref{TH:local-lq}, (ii). The main statement is
part (ii), and (i) serves to highlight a difference between
$L^q$ and $\FT L^q$ for such problems. The proof of (ii)
will rely on some properties of Fourier integral operators
in an essential way.

Let us now discuss an implication of the boundedness result
for the regularity properties of Fourier integral operators.
We note that since $\psi^*$ is bounded on $L^q(\Rn)$ in the
assumptions of Theorem \ref{TH:local-lq}, it follows that
$I_\psi$ is bounded on $\FT L^q(\Rn)$. By an argument
similar to the one that we will give in the proof of Theorem
\ref{TH:local-lq} this implies that $I_\psi$
is continuous from $(\FT L^q)_{comp}$ to $(\FT L^q)_{loc}$,
so that Theorem \ref{TH:local-lq} also follows
if we use the equalities
from Theorem \ref{TH:mw}. This is related to the
question of the local boundedness of
Fourier integral operators on $\FT L^q$.
Let $T$ be defined by
$$
 Tf(x)=\int_\Rn e^{i\Phi(x,\xi)} a(x,\xi) \widehat{f}(\xi)
 d\xi,
$$
where $\Phi$ is a non-degenerate real-valued phase function.
In \cite{CNR08}, it was shown that if
the phase function is non-degenerate and
homogeneous of order one and if
the amplitude $a(x,\xi)$ is compactly supported in $x$ and
belongs to the symbol class $S^m_{1,0}$ with $m\leq -n|1/q-1/2|$,
then $T$ is bounded on
$(\FT L^q)_{comp}$.
%First, we observe that this result means
%that operators $T$ of the form
%$$
% Tf(x)=\int_\Rn\int_\Rn e^{i(\Phi(x,\eta)-y\cdot\eta)}
% b(y,\eta) \widehat{f}(y) dy  d\eta,
%$$
%where $b\in S^m_{1,0}$ with $m\leq -n|1/p-1/2|$, are
%locally continuous on $\FT L^q$, i.e.
%continuous from $(\FT L^q)_{comp}$ to $(\FT L^q)_{loc}$.
Moreover, they showed the order $m$ to be sharp
for a special choice of the phase function $\Phi(x,\xi)$.
However,
Theorem \ref{TH:local-lq} implies that if we take
the amplitude $a$ in the class $S^0_{0,0}$, and
the phase function corresponding to the canonical transforms,
operator $T$ is still locally continuous on $\FT L^q$, i.e.
continuous from $(\FT L^q)_{comp}$ to $(\FT L^q)_{loc}$.
We note the inclusion
$S^0_{0,0}(\Rn)\subset M^{\infty,1}(\R^{2n})$
(here always $S^0_{0,0}(\Rn\times\Rn)$ is defined as the set of all
smooth $a\in C^\infty(\R^{2n})$ such that
$|\partial_x^\alpha\partial_\xi^\beta a(x,\xi)|\leq
C_{\alpha\beta}$ for all multi-indices $\alpha,\beta$ and all
$x,\xi\in\Rn$). Thus, we have the following result:

\begin{thm}\label{cor:FIOs}
Let $\psi : \R^n \to \R^n$ be such that
$\psi^*$ is bounded on $L^q(\Rn)$ and let
$a\in M^{\infty,1}(\R^{2n})$.
Then the operator
$$
 Tf(x)=\int_\Rn\int_\Rn e^{i(x\cdot\xi-y\cdot\psi(\xi))}
 a(x,\xi) {f}(y) dy  d\xi
$$
is locally continuous on $\FT L^q$, $M^{p,q}$ and
$W^{p,q}$, for all $1\leq p,q<\infty$.
\end{thm}

We note that this result is also true for $p=\infty$ and $q=\infty$
if we use the modification as in Remark \ref{modified-def}.

\par

We will also discuss non-affine
transforms which induce the globally bounded changes of
variables on $M^{p,q}$ or $W^{p,q}$.
Note that such transforms must not be a $C^1$-mappings in view of
Theorem \ref{TH:B-H}.
Moreover, we will show that the Beurling--Helson type theorem
fails if we allow derivatives of $\psi$ to have singularities
of types important for applications to partial differential
equations.
One example of this is Theorem \ref{thm:hom}.
In fact, to prove the conclusion of part (ii) of Theorem
\ref{thm:hom} we will use the sharpness results on the
$L^q$ boundedness of Fourier integral operators established
in \cite{Ru99}.

Finally, we establish several positive results for homogeneous
changes of variable which may have more singularities than
only at the origin.
We investigate properties for mappings of the form
\begin{equation}\label{mapping1}
f(x)\mapsto f(S(x)+T(|x_1|,\dots |x_n|)),
\end{equation}
when acting on modulation spaces or Wiener amalgam spaces.
Here $S$ and $T$ are linear mappings on $\rr n$ such that
\begin{equation}\label{mapping1A}
x\mapsto S(x)+T((-1)^{j_1}x_1,\dots ,(-1)^{j_n}x_n)
\end{equation}
is a bijection on $\rr n$, for each choice of $j_1,\dots j_n\in \{
0,1\}$. In particular, the following situations are covered by
\eqref{mapping1}:
\begin{enumerate}
\item[(i)]  $f(x)\mapsto f(|x_1|,\dots ,|x_n|)$, which follows by
choosing
$$
S=0\quad \text{and} \quad T=\operatorname{Id} _{\rr n};
$$

\vrum

\item[(ii)] $f(x)\mapsto f(x_1,\dots ,x_{n-1},|x_n|)$,
which follows by choosing
$$
S(x)=(x_1,\dots ,x_{n-1},0)\quad \text{and} \quad
T=\operatorname{Id} _{\rr n} -S;
$$

\vrum

\item[(iii)] $f(x)\mapsto f(x_1,\dots ,x_{n-1},|x_1|+\cdots +|x_n|)$,
which follows by choosing
$$
S(x)=(x_1,\dots ,x_{n-1},0)\quad \text{and} \quad
T(x)=(0,\dots ,0,x_1+\cdots +x_n).
$$
\end{enumerate}

\par

For such mappings we have the following result:

\par

\begin{thm}\label{homogen}
Assume that $p,q\in (1,\infty )$, and assume that
$S$ and $T$ are linear mappings on $\rr n$ such that
for each $j_1,\dots j_n\in \{ 0,1\}$, the map \eqref{mapping1A} is
bijective. Then the map \eqref{mapping1} from $\mathscr S(\rr n)$ to
$\mathscr S'(\rr n)$ extends uniquely to continuous mappings on
$M^{p,q}(\rr n)$ and on $W^{p,q}(\rr n)$.
\end{thm}

Theorem \ref{homogen} says in particular that
if a homogeneous of order one function $\psi$ has more
singularities than only at the origin, then $\psi^*$ may
still be bounded on modulation and Wiener amalgam spaces
$M^{p,q}$ and $W^{p,q}$. We note that this type of statement
on $\FT L^q$ appeared in \cite{LO94} while the case of
$M^{p,p}$ was analysed in \cite{Ok08}.

All these results will be proved in Section 4.
Some results of this paper
were partially announced by authors in \cite{RSTT}.

%%%==================================================================
%%%==================================================================
\section{Preliminaries}
Let $\calS(\R^n)$ and $\calS'(\R^n)$ be the Schwartz spaces of
all rapidly decreasing smooth functions
and tempered distributions,
respectively.
We define the Fourier transform $\FT f$
and the inverse Fourier transform $\FT^{-1}f$
of $f \in \calS(\R^n)$ by
\[
\FT f(\xi)
=\widehat{f}(\xi)
=\int_{\R^n}e^{-ix \cdot \xi}\, f(x)\, dx
\quad \text{and} \quad
\FT^{-1}f(x)
=\frac{1}{(2\pi)^n}
\int_{\R^n}e^{ix\cdot \xi}\, f(\xi)\, d\xi.
\]

\par

We introduce modulation spaces
based on Gr\"ochenig in \cite{Grochenig}.
Fix a function $\varphi \in \calS(\R^n)\setminus 0$
(called the \emph{window function}).
Then the short-time Fourier transform $V_{\varphi}f$ of
$f \in \calS'(\R^n)$ with respect to $\varphi$
is defined by
\[
V_{\varphi}f(x,\xi)
=\langle f, M_{\xi}T_x \varphi\rangle
=\int_{\R^n}f(t)\, \overline{\varphi(t-x)}\, e^{-i\xi\cdot t}\, dt
\]
for $x, \xi \in \R^n$,
where $M_{\xi}\varphi(t)=e^{i\xi \cdot t}\varphi(t)$
and $T_x \varphi(t)=\varphi(t-x)$.
%and $\langle\cdot,\cdot\rangle$ denotes the inner product on $L^2(\R^n)$.
We note that,
for $f \in \calS'(\R^n)$,
$V_{\varphi}f$ is continuous on $\R^{2n}$
and $|V_{\varphi}f(x,\xi)|\le C(1+|x|+|\xi|)^N$
for some constants $C,N \ge 0$
(\cite[Theorem 11.2.3]{Grochenig}).

\par

Let $1\le p,q \le \infty$. Then we let $L^{p,q}_1(\rr {2n})$ be the
set of all $F\in L^1_{loc}(\rr {2n})$ such that $\nm
F{L^{p,q}_1}<\infty$, where
\begin{alignat*}{2}
&\|F\|_{L^{p,q}_1}
=\left\{ \int_{\R^n} \left(
\int_{\R^n} |F(x,\xi)|^{p}\, dx
\right)^{q/p} d\xi \right\}^{1/q},&
\qquad 1 &\le p,q<\infty,
\\[1ex]
&\|F\|_{L^{\infty,q}_1}
=\left\{\int_{\R^n}
\left(\underset {x \in \R^n} {\operatorname{ess\, sup}}|F(x,\xi)|
\right)^q d\xi \right\}^{1/q},&
\qquad 1 &\le q<\infty,
\\[1ex]
&\|F\|_{L^{p,\infty}_1}
=\underset {x \in \R^n} {\operatorname{ess\, sup}}\left(\int_{\R^n}
|F(x,\xi)|^p \, dx \right)^{1/p},&
\qquad 1 &\le p<\infty,
\\[1ex]
&\|F\|_{L^{\infty,\infty}_1}
=\underset {x,\xi \in \R^n}{\operatorname{ess\, sup}}
|F(x,\xi)|.& &
\end{alignat*}
The modulation space $M^{p,q}(\R^n)$ consists of all $f \in
\calS'(\R^n)$ such that $V_\fy f(x,\xi )\in L^{p,q}_1(\rr {2n})$, i.{\,}e.
$M^{p,q}(\R^n)$ consists of all $f \in \calS'(\R^n)$ such that
$\|f\|_{M^{p,q}} \equiv \nm {V_\varphi f}{L^{p,q}_1}$ is finite.
If $p=q$, we simply write $M^p$ instead of $M^{p,p}$.
We note that $M^{2,2}(\R^n)=L^2(\R^n)$,
$M^{p,q}(\R^n)$ is a Banach space under the norm $\|\cdot
\|_{M^{p,q}}$, $\calS(\R^n)$ is dense in $M^{p,q}(\R^n)$ if $1 \le
p,q<\infty$, and $M^{p_1,q_1}(\R^n) \hookrightarrow M^{p_2,q_2}(\R^n)$
if $p_1 \le p_2$ and $q_1 \le q_2$
(cf. Propositions 11.3.1, 11.3.4, 11.3.5 and Theorem 12.2.2 in
\cite{Grochenig}).
The definition of $M^{p,q}(\R^n)$ is independent
of the choice of the window function
$\varphi \in \calS(\R^n)\setminus 0$,
that is, different window functions
yield equivalent norms
(\cite[Proposition 11.3.2]{Grochenig}).
We denote by $p' \in [1,\infty]$ the conjugate exponent
of $p \in [1,\infty]$, i.{\,}e. $1/p+1/p'=1$.

\begin{rem}[{\cite[Lemma 2.2]{BGHO}, \cite[Lemma 3.2]{To04}}]
\label{modified-def}
Let $1 \le p,q \le \infty$,
and let $\calM^{p,q}(\R^n)$ be the completion of $\calS(\R^n)$
under the norm $\|\cdot\|_{M^{p,q}}$.
Then the following are true:
\begin{enumerate}
\item[(i)]
if $1 \le p,q<\infty$, then $\calM^{p,q}=M^{p,q}$;
\item[(ii)]
%if $1 \le p,q<\infty$
%then $(\calM^{\infty,q})'=\calM^{1,q'}$
%and $(\calM^{p,\infty})'=\calM^{p',1}$,
%and $(\calM^{\infty,\infty})'=\calM^{1,1}$.
%
%  Are these identities correct? Here below I have reformulate
%  (ii) what I believe is correct.
%
if $1 \le p,q<\infty$
then $(\calM^{\infty,q})'=M^{1,q'}$
and $(\calM^{p,\infty})'=M^{p',1}$,
and $(\calM^{\infty,\infty})'=M^{1,1}$.
\end{enumerate}
\end{rem}

Next, we discuss Wiener amalgam spaces. We let $L^{p,q}_2(\rr {2n})$
be the set of all $F\in L^1_{loc}(\rr {2n})$ such that $\nm
F{L^{p,q}_2}<\infty$, where
\begin{alignat*}{2}
&\|F\|_{L^{p,q}_2}
=\left \{ \int_{\R^n} \left(
\int_{\R^n} |F(x,\xi)|^{q}\, d\xi
\right ) ^{p/q}\,  dx \right\}^{1/p}, &
\qquad 1 &\le p,q<\infty,
\\[1ex]
&\|F\|_{L^{\infty,q}_2}
=\underset {x \in \R^n}{\operatorname{ess\, sup}}\left( \int_{\R^n}
|F(x,\xi)|^q \, d\xi \right)^{1/q},&
\qquad 1 &\le q<\infty,
\\[1ex]
&\|F\|_{L^{p,\infty}_2}
=\left\{ \int_{\R^n} \left(
\underset {\xi \in \R^n}{\operatorname{ess\,
sup}}|F(x,\xi)|\right)^p dx \right \} ^{1/p},&
\qquad 1 &\le p<\infty,
\\[1ex]
&\|F\| _{L^{\infty ,\infty}_2}
=\underset {x,\xi \in \R^n}{\operatorname{ess\, sup}}
|F(x,\xi )|.&  &
\end{alignat*}
Obviously, if $F(x,\xi )=G(\xi ,x)$, then $F\in L^{p,q}_1$ if and only
if $G\in L^{q,p}_2$. We set $\nm f{W^{p,q}}=\nm {V_\fy f}{L^{p,q}_2}$. The 
Wiener
amalgam space $W^{p,q}(\R^n)$
consists of all $f \in \calS'(\R^n)$ such that
$\|f\|_{W^{p,q}}<\infty$. (Note that the general definition of Wiener
amalgam spaces in \cite{Fe2} permits function and distribution spaces
which are not considered here.) Since
\begin{equation}\label{fourtransfrel}
\left|V_{\varphi}f(x,\xi)\right|
=(2\pi)^{-n}\left|V_{\widehat{\varphi}}\widehat{f}(\xi,-x)\right|,
\end{equation}
we see that
\begin{equation}\label{(2.1)}
 \|f\|_{W^{p,q}}\asymp\|\widehat{f}\|_{M^{q,p}}.
\end{equation}
This implies that
the definition of $W^{p,q}(\R^n)$ is independent
of the choice of the window function
$\varphi \in \calS(\R^n)\setminus 0$,
since the modulation space $M^{q,p}(\R^n)$ is so.
By the same reason, we also have
$W^{p_1,q_1}(\R^n) \hookrightarrow W^{p_2,q_2}(\R^n)$
if $p_1 \le p_2$ and $q_1 \le q_2$, and other properties
similar to those of $M^{p,q}$.

\par

In Appendix A, we consider general modulation spaces and weighted
versions of Wiener amalgam spaces.

\par

%%%==================================================================
%%%==================================================================
\section{Proofs of the main results, and some further remarks}
%%%==================================================================

In this section we prove our results. The following proposition is
needed in the proof of Theorem \ref{TH:mw}.

\begin{prop}\label{3.1}
Let $1 \le p,q \le \infty$ and $\Omega$ be a compact subset of $\R^n$.
Then the following are true:
\begin{enumerate}
\item[(i)]
there exists a constant $C>0$ such that
\[
\|f\|_{M^{p,q}} \le
C|\widetilde{\Omega}|^{\max(0,1/q-1/p)}\|f\|_{W^{p,q}}
\quad \text{for all $f \in W^{p,q}(\R^n)$
with $\mathrm{supp}\, f \subset \Omega$;}
\]
\item[(ii)]
there exists a constant $C>0$ such that
\[
\|f\|_{W^{p,q}} \le
C|\widetilde{\Omega}|^{\max(0,1/p-1/q)}\|f\|_{M^{p,q}}
\quad \text{for all $f \in M^{p,q}(\R^n)$
with $\mathrm{supp}\, f \subset \Omega$,}
\]
\end{enumerate}
where $\widetilde{\Omega}=\{x \in \R^n : \mathrm{dist}(x,\Omega)<1\}$,
and $C>0$ is independent of $\Omega$.
\end{prop}
%%%==================================================================
\begin{proof}
Let $\Omega$ be a compact subset of $\R^n$,
and let $\varphi \in \calS(\R^n)\setminus 0$
with $\mathrm{supp}\, \varphi \subset B(0,1)$,
where $B(0,1)$ is the open ball with radius $1$
centred at the origin.
Assume that $f \in \calS'(\R^n)$
with $\mathrm{supp}\, f \subset \Omega$. Then,
$$
\mathrm{supp}\, V_{\varphi}f(\cdot,\xi) \subset
\widetilde{\Omega}=\Omega+B(0,1) =\{x \in \R^n :
\mathrm{dist}(x,\Omega)<1\}
$$
for all $\xi \in \R^n$.

\par

We first consider the case $p \le q$.
By Minkowski's inequality,
\begin{align*}
\|f\|_{M^{p,q}}
&=\left\{ \int_{\R^n} \left(
\int_{\R^n} |V_{\varphi}f(x,\xi)|^p\, dx\right)^{q/p} d\xi \right\}^{1/q}
\\
&\le \left\{ \int_{\R^n} \left(
\int_{\R^n} |V_{\varphi}f(x,\xi)|^q\, d\xi\right)^{p/q} dx \right\}^{1/p}
=\|f\|_{W^{p,q}}.
\end{align*}
On the other hand, by H{\"o}lder's inequality,
\begin{align*}
&\|f\|_{W^{p,q}}
=\left\{ \int_{\widetilde{\Omega}} \left(
\int_{\R^n} |V_{\varphi}f(x,\xi)|^q\, d\xi\right)^{p/q}
dx \right\}^{1/p}
\\
&\le \left[ \left\{\int_{\widetilde{\Omega}} \left(
\int_{\R^n} |V_{\varphi}f(x,\xi)|^q\, d\xi\right) dx
\right\}^{p/q}|\widetilde{\Omega}|^{1-p/q} \right]^{1/p}
=|\widetilde{\Omega}|^{1/p-1/q}\|f\|_{M^{q,q}}.
\end{align*}
Since $M^{p,q} \hookrightarrow M^{q,q}$,
we see that
\[
\|f\|_{W^{p,q}}
\le |\widetilde{\Omega}|^{1/p-1/q}\|f\|_{M^{q,q}}
\le C|\widetilde{\Omega}|^{1/p-1/q}\|f\|_{M^{p,q}}.
\]

\par

We next consider the case $p \ge q$.
In the same way as in the case $p \le q$,
by Minkowski's inequality,
we have $\|f\|_{W^{p,q}} \le \|f\|_{M^{p,q}}$.
On the other hand,
by H\"older's inequality,
we see that
\begin{align*}
&\|f\|_{M^{q,q}}
=\|f\|_{W^{q,q}}
=\left\{ \int_{\widetilde{\Omega}} \left(
\int_{\R^n} |V_{\varphi}f(x,\xi)|^q\, d\xi\right)
dx \right\}^{1/q}
\\
&\le \left[ \left\{\int_{\widetilde{\Omega}} \left(
\int_{\R^n} |V_{\varphi}f(x,\xi)|^q\, d\xi\right)^{p/q}
dx \right\}^{q/p}|\widetilde{\Omega}|^{1-q/p}\right]^{1/q}
=|\widetilde{\Omega}|^{1/q-1/p}\|f\|_{W^{p,q}}.
\end{align*}
Hence, it follows from the embedding
$M^{q,q} \hookrightarrow M^{p,q}$ that
\[
\|f\|_{M^{p,q}}\le C\|f\|_{M^{q,q}}
\le C|\widetilde{\Omega}|^{1/q-1/p}\|f\|_{W^{p,q}}.
\]
The proof is complete.
\end{proof}
%%%==================================================================
%We are now ready to prove Theorems \ref{TH:mw}--\ref{TH:local-lq}.
%\begin{proof}[Proof of Theorem \ref{TH:mw} (Ruzhansky).]
%The proof of \eqref{EQ:mw-eq} is simple if we use the following
%expressions for norms in modulation and Wiener amalgam spaces:
%\begin{equation}\label{EQ:alt-norms}
%||f||_{M^{p,q}}=\left|\left| ||\Phi(D-k)u||_{L^p_x}
%\right|\right|_{l^q_k},\quad
%||f||_{W^{p,q}}=\left|\left| ||\Phi(D-k)u||_{l^q_k}
%\right|\right|_{L^p_x},
%\end{equation}
%for some non-zero $\Phi\in C_0^\infty(\Rn)$, and where we
%can take the support of $\Phi$ to be contained in the open ball
%with radius $1/2$ centred at the origin.
%Now, we can observe that if $\Omega\in\Rn$ is a compact
%set contained in a ball with radius $1/2$
%(but centred at any point) and
%if $\widehat{f}\in\calE'(\Omega)$ then
%we get equalities
%$||f||_{M^{p,q}}=||\Phi(D-k_0)u||_{L^p_x}$ and
%$||f||_{W^{p,q}}=||\Phi(D-k_0)u||_{L^p_x}$, for some
%$k_0$. Moreover, if the support of $\widehat{f}$ is
%an arbitrary compact set, we get finite sums of these
%expressions, implying the second line in \eqref{EQ:mw-eq}.
%The first line follows from the second one by taking the
%Fourier transform. Finally, estimates in \eqref{EQ:mw-norms}
%follow from Proposition \ref{3.1}.
%\end{proof}
%%%==================================================================

We are now ready to prove Theorems \ref{TH:mw}--\ref{TH:local-lq}.

\begin{proof}[Proof of Theorem \ref{TH:mw}.]
The proof of \eqref{EQ:mw-eq} is simple if we use the following
expressions for norms in modulation and Wiener amalgam spaces:
\begin{equation}\label{EQ:alt-norms}
\|f\|_{M^{p,q}}\asymp\left\| \|\Phi(D-k)u\|_{L^p_x}
\right\|_{l^q_k},\quad
\|f\|_{W^{p,q}}\asymp\left\| \|\Phi(D-k)u\|_{l^q_k}
\right\|_{L^p_x},
\end{equation}
%for some non-zero $\Phi\in C_0^\infty(\Rn)$, and where we
%can take the support of $\Phi$ to be contained in the open cube $(-1,1)^n$.
where $\Phi \in C_0^\infty(\R^n)$ satisfies
$\mathrm{supp}\, \Phi \subset [-1,1]^n$ and
$\sum_{k \in \Z^n}\Phi(\xi-k)\equiv 1$.
Now, we can observe that if $\Omega\in\Rn$ is a compact
set contained in an open cube with side-length $2$
(but centred at any point) and
if $\widehat{f}\in\calE'(\Omega)$ then
we get
\[
\|f\|_{M^{p,q}}\asymp\sum_{\substack{|k_j-k_{0,j}|\le 2 \\ j=1,\dots,n}}
\|\Phi(D-k_0)u\|_{L^p}, \quad
\|f\|_{W^{p,q}}\asymp\sum_{\substack{|k_j-k_{0,j}|\le 2 \\ j=1,\dots,n}}
\|\Phi(D-k_0)u\|_{L^p}
\]
for some $k_0$, where $k=(k_1,\dots,k_n) \in \Z^n$
and $k_0=(k_{0,1},\dots,k_{0,n}) \in \Z^n$.
Moreover, if the support of $\widehat{f}$ is
an arbitrary compact set, we get finite sums of these
expressions, implying the second line in \eqref{EQ:mw-eq}.
The first line follows from the second one by taking the
Fourier transform. Finally, estimates in \eqref{EQ:mw-norms}
follow from Proposition \ref{3.1}.
\end{proof}
%%%==================================================================

\par

We note that the identity
$$
M^{p,q}\cap \calE'=\FT L^q\cap \calE'
$$
in \eqref{EQ:mw-eq} also follows from Remark 1.3 (4) in \cite{CT1},
and it was announced in several conferences by the authors, including
``Mathematical modeling of wave phenomena 05'', in V{\"a}xj{\"o}, Sweden
(see also Remark \ref{thmagain} below). A more recent 
alternative proof of the
latter equality can be found in \cite[Lemma 1]{Ok08}. In this context
we pay attention to the simplicity of the proof of Theorem \ref{TH:mw}
here as above, based on our choice of using the norm \eqref{EQ:alt-norms}
instead of short time Fourier transforms. The equality
$$
W^{p,q}\cap \calE'=\FT L^q\cap \calE'
$$
concerning Wiener amalgam spaces $W^{p,q}$ appears to be new.

\par

In Remark \ref{thmagain} below we give an extension of
\eqref{EQ:mw-eq}, based on a different technique compared to the
proof of Theorem \ref{TH:mw}, and which involves weighted
spaces. These considerations are dependent on some multiplication and
convolution properties for modulation spaces which we shall discuss
now.

\par

\begin{rem}\label{thmagain}
In addition to two proofs contained in this paper,
there are also other ways to obtain the inclusion \eqref{EQ:mw-eq}. In
fact, we can use the multiplication properties for modulation
spaces in Appendix A to obtain the latter inclusion in a more general
context involving spaces of the form $M^{p,q}_{(\omega )}$ and
$W^{p,q}_{(\omega )}$, which are now dependent of the weight function
$\omega \in \mathscr P(\rr {2n})$ (cf. Appendix A for precise
definitions of $\mathscr P(\rr {2n})$, $M^{p,q}_{(\omega
)}$ and $W^{p,q}_{(\omega )}$.) We claim that
\begin{equation}\tag*{(\ref{EQ:mw-eq})$'$}
\begin{alignedat}{2}
&M^{p,q}_{(\omega )}\cap\calE'=W^{p,q}_{(\omega )}\cap\calE'=\FT
L^q_{(\omega _0)}\cap\calE',\quad \omega (x,\xi )&= \omega _0(\xi )\in
\mathscr P(\rr n),
\\
&M^{p,q}_{(\omega )}\cap \FT\calE'=W^{p,q}_{(\omega )}\cap
\FT\calE'=L^p_{(\omega )}\cap \FT\calE',\omega (x,\xi )&= \omega
(x)\in \mathscr P_0(\rr n),
\end{alignedat}
\end{equation}
with equivalence of norms.

\par

Indeed, assume that $f\in M^{\infty ,q}_{(\omega )}\cap \calE'$,
$\chi \in C_0^\infty$ is equal to $1$ in the support of $f$, and
$v_0\in \mathscr P(\rr {n})$ is such that $\omega _0$ is
$v_0$-moderate. Then $\chi \in M^{1,1}_{(v)}$, where $v(x,\xi
)=v_0(\xi )$. Hence Proposition A.1 gives
$$
f =f\, \chi \in M^{\infty ,q}_{(\omega )}\cdot M^{1,1}_{v}\subseteq
M^{1,q}_{(\omega )}.
$$
This proves that $M^{p,q}_{(\omega )}\cap \calE'$ is independent of
$p$. By similar arguments it follows that $W^{p,q}_{(\omega )}\cap
\calE'$ is independent of $p$. The first two equalities in
\eqref{EQ:mw-eq}$'$ is now a consequence of {\rm{(A.1)}}.

\par

The last part of \eqref{EQ:mw-eq}$'$ follows by similar arguments,
using Proposition A.2 instead of Proposition A.1. Alternatively, the
second line in \eqref{EQ:mw-eq}$'$ follows from the first line and
the Fourier inversion formula.
\end{rem}

\par

Before proving Theorem \ref{TH:B-H}, let us point out the following
immediate consequence of Proposition A.1 in Appendix, which we will
use to investigate localisation properties.

\par

\begin{prop}\label{prop:mult}
Assume that $p_j ,q_j\in [1, \infty]$ for $j=0,1,2$ satisfy
$$
\frac 1{p_1}+\frac 1{p_2}=\frac 1{p_0}\quad \text{and}\quad \frac
1{q_1}+\frac 1{q_2}=1+\frac 1{q_0}.
$$
Then the map $(f_1,f_2)\mapsto f_1\cdot f_2$ on $\mathscr S(\rr n)$
extends to continuous mappings from $M^{p_1,q_1}(\rr n)\times
M^{p_2,q_2}(\rr n)$ to $M^{p_0,q_0}(\rr n)$, and from $W^{p_1,q_1}(\rr
n)\times W^{p_2,q_2}(\rr n)$ to $W^{p_0,q_0}(\rr n)$. Furthermore,
each modulation space or Wiener amalgam space is an $M^{\infty
,1}$-module under multiplication.
\end{prop}

\par

\begin{proof}
The asserted mapping property follows from Proposition A.1, or from
Theorem 3 in \cite{Fe3} for modulation spaces and from Theorem 2.4 in
\cite{To04} for Wiener amalgam spaces. By letting $p_1=\infty$ and
$q_1=1$, it follows that $p_2=p_0$ and $q_2=q_0$. Hence $M^{p_2,q_2}$
is an $M^{\infty ,1}$ module and $W^{p_2,q_2}$ is an $W^{\infty ,1}$
module. The asserted module properties now follows from these
relations and the fact that $M^{\infty ,1}\subseteq W^{\infty ,1}$.
\end{proof}

In what follows we let $\mathscr M_\uppsi$ denote the multiplication
operator $\mathscr M_\uppsi f=\uppsi \cdot f$, for appropriate
functions or distributions $f$ and $\uppsi$.

\begin{prop}\label{op-equiv}
Let $1 \le p,q < \infty$
and $\uppsi \in \calS'(\Rn)$.
Then the following is true:
\begin{itemize}
\item[(i)]
$\MT_{\uppsi}$ is bounded on $M^{p,q}(\Rn)$
if and only if $\uppsi (D)$ is bounded on $W^{q,p}(\Rn)$;
\item[(ii)]
$\uppsi (D)$ is bounded on $M^{p,q}(\Rn)$
if and only if $\MT_{\uppsi}$ is bounded on $W^{q,p}(\Rn)$.
\end{itemize}
\end{prop}

\begin{proof}
Assume that $\uppsi (D)$ is bounded on $W^{q,p}(\Rn)$.
By \eqref{(2.1)},
we see that
\[
\|\MT_{\uppsi}f\|_{M^{p,q}} \asymp
\|\FT^{-1}[\MT_{\uppsi}f]\|_{W^{q,p}}=\|\uppsi (D)[\FT^{-1}f]\|_{W^{q,p}}.
\]
Hence,
\begin{align*}
\|\MT_{\uppsi}f\|_{M^{p,q}}
&\le C\|\uppsi (D)[\FT^{-1}f]\|_{W^{q,p}}
\\
&\le C\|\uppsi (D)\|_{\calL(W^{q,p})}\|\FT^{-1}f\|_{W^{q,p}}
\le C\|\uppsi (D)\|_{\calL(W^{q,p})}\|f\|_{M^{p,q}}.
\end{align*}
In the same way,
we can prove the others in Proposition \ref{op-equiv}.
\end{proof}

\begin{rem}
Propositions \ref{prop:mult} and \ref{op-equiv}
with $p=\infty$ or $q=\infty$ hold under the modification
as in Remark \ref{modified-def}.
\end{rem}

%\begin{rem}\label{Rem:norm-chi}
%We note that if the characteristic function
%$\chi_\Omega$ satisfies $\chi_\Omega\in M^{\infty,1}(\Rn)$,
%then the multiplication by $\chi_\Omega$ is bounded on
%$M^{p,q}(\Rn)$ and $W^{p,q}(\Rn)$ by
%Proposition \ref{prop:mult}, so that inequalities
%\eqref{EQ:mw-norms} imply that the operator norm
%of multiplication by $\chi_\Omega$ as a mapping
%from $M^{p,q}(\Rn)$ to $W^{p,q}(\Rn)$, and from
%$W^{p,q}(\Rn)$ to $M^{p,q}(\Rn)$, is bounded by
%$C(1+|\Omega|)^{1-1/p}||\chi_\Omega||_{M^{\infty,1}(\Rn)}.$
%
%{\bf by the way, is it well-known when
%$\chi_\Omega\in M^{\infty,1}(\Rn)$? (if the boundary
%of $\Omega$ is, say, smooth) or, also, what are the
%trace theorems in modulation spaces - any loss of regularity
%like for Sobolev spaces? - I'm just curious!}
%\end{rem}

\begin{rem}\label{Rem:norm-chi}
We note that $\mathscr M_\uppsi $ in Proposition \ref{op-equiv} is
bounded on $M^{p,q}$ and on $W^{p,q}$ when 
$\uppsi \in M^{\infty ,1}$
by Proposition \ref{op-equiv}.
In this context we note that if $\uppsi$ is a
characteristic function for sets with non-zero Lebesgue measure, then
$\uppsi \notin M^{\infty ,1}$ and  $\uppsi \notin W^{\infty ,1}$,
since $M^{\infty ,1}\subseteq W^{\infty ,1}$ are contained in the set
of continuous functions (in the distributional sense, see
Remark \ref{Rem:W-cont} for the proof).

\par

On the other hand, $\chi_{[-1,1]^n}(D)$ is bounded on $L^p(\R^n)$ and
on $M^{p,q}$ with $1<p<\infty$,
but $\chi_{[-1,1]^n} \not\in M^{\infty,1}(\R^n)$
(cf. \cite[Proposition 3.6]{Duo} and \cite{BGGO05}).
\end{rem}

\begin{rem}\label{Rem:W-cont}
We note that $W^{\infty ,1}$ is contained in the set
of continuous functions (in the distributional sense).
Since $M^{\infty ,1}\subseteq W^{\infty ,1}$ the same
is automatically true for the space $M^{\infty ,1}$.

Since we were not able to find this statement in the
literature, we will now give a simple justifying argument.
Assume that $f\in W^{\infty ,1}$ and choose a window 
function $\chi\in C^\infty$ such that $\chi(0)=1$. Then function
$$\xi\mapsto \FT (f \chi ( \cdot  -x_0))( \xi )$$
belongs to $L^1$ for every fixed $x_0$. Hence its 
inverse Fourier transform 
$$
y   \mapsto  \FT^{-1} (\FT (f \chi ( \cdot  -x_0)))(y)=
f(y) \chi (y-x_0)
$$
is continuous. Since $\chi (y-x_0)$ is smooth and non-zero 
around $y=x_0$, it follows that 
$f(y)$  is continuous around $x_0$. Since 
$x_0$ was arbitrarily chosen it follows that 
$f$ is continuous everywhere.
\end{rem}

\begin{proof}[Proof of Theorem \ref{TH:B-H}.]
Let us first consider condition (ii) in which case
the statement is a straightforward
consequence of the same
statement for $\FT L^q$ (Beurling--Helson type theorem;
see also \cite{BH, LO94, Se76} for the case of $\FT L^q$ and
\cite{Ok08} for the case of $M^{p,q}$).
The case of $W^{p,q}$ in (ii) follows if we use the equality
\eqref{EQ:mw-eq} for the localised versions.

Now, since conditions (i) and (iii) as well as (ii) and (iv)
are equivalent, respectively, in view of relation
\eqref{EQ:ch-can}, it is enough to show that assumption (i)
implies (ii). But this follows immediately from Proposition
\ref{prop:mult} and the inclusion
$C_0^\infty\subset M^{\infty,1}.$
\end{proof}

\begin{proof}[Proof of Theorem \ref{TH:local-lq}.]
Here, because (i) and (ii) are equivalent in view of
\eqref{EQ:ch-can},
all we need to show is the boundedness
of $\chi_1(x)I_\psi \chi_2(x)$ on $M^{p,q}$,
where $\chi_1,\chi_2\in C^\infty_0(\Rn)$.
Equivalently, we may show the boundedness of
$\chi_1(D)\psi^* \chi_2(D)$ on $W^{p,q}$
if we use the relation \eqref{(2.1)}. The latter is
induced by
the $L^q(\Rn)$--boundedness of $\chi_1(D)$ and $\chi_2(D)$
($1\le q \le\infty$) due to Young's inequality
(note that the kernels are in $L^1$).
\end{proof}

\begin{proof}[Proof of Theorem \ref{cor:FIOs}.]
Indeed, since pseudo-differential operator $a(x,D)$ with
symbol $a\in M^{\infty,1}$ is bounded on the modulation
space $M^{p,q}(\Rn)$, it is also locally continuous on
$M^{p,q}(\Rn)$ in view of Proposition \ref{prop:mult}.
But then $T=a(x,D)\circ I_\psi$ is locally continuous
on $M^{p,q}(\Rn)$ by Theorem \ref{TH:local-lq}, and hence
also on $\FT L^q$ and $W^{p,q}$ by Theorem \ref{TH:mw}.
\end{proof}

\begin{proof}[Proof of Theorem \ref{thm:hom}.]
(i) It follows that
$\psi^{-1}\in C^\infty(\Rn\backslash 0)$
is also positively homogeneous
or order one. Hence its derivative $D\psi^{-1}$
is homogeneous of order zero and hence bounded on $\Rn$.
Consequently, $\psi^*$ is bounded on $L^q(\Rn)$ for all
$1\leq q\leq\infty$ and statement (i) follows from Theorem
 \ref{TH:local-lq}.

(ii) Suppose now that $\psi^*$ is continuous on
$\calL(\FT L^q(\Rn))$ for
$1\leq q\leq \infty$, $q\not=2$.
Then it follows that $I_\psi$ is bounded (and hence also
locally bounded) on $L^q(\Rn)$
(the conclusion that $I_\psi$ is locally bounded on $L^q(\Rn)$
is also true if 
$\psi^*$ is Fourier-locally continuous on
$\calL(\FT L^q(\Rn))$).
In turn, this implies that
$\psi$ must be a linear function, if we use the critical
orders for the $L^q$--boundedness of Fourier integral operators
obtained in \cite{Ru99}. For completeness (and also to include
boundary cases $q=1$ and $q=\infty$), let us give this argument
in more detail. We can assume that $1\leq q<2$ since the case
$2<q\leq\infty$ follows by considering the adjoints.

Let $\psi:\Rn\to\Rn$ and assume that
$\psi\in C^\infty(\Rn\backslash 0)$ and that it
is positively homogeneous of
order one.
Let us define $k$ by setting
$$\max_{\xi\in\Rn\backslash 0, \;y\in\Rn}
\rank \nabla^2\psi(\xi)y=n-k.$$
First we observe that the canonical relation
of the Fourier integral
operator $I_\psi$
is given by
$$\Lambda=\b{\p{\nabla\psi(\xi)y,\xi,y,\psi(\xi)}: y\in\Rn,
 \xi\in\Rn\backslash 0}.$$
Its projection $\Sigma=\pi(\Lambda)$ to the base space
$\Rn\times\Rn$ is a set of dimension $\leq 2n-k$.
At points where $\rank \nabla^2\psi(\xi)y=n-k$ this is a smooth
manifold of dimension $2n-k$ (with its conormal
bundle equal to $\Lambda$). Let it be locally given by the
set of defining equations $h_j(x,y)=0$, $j=1,\ldots,k$, with
$\nabla h_1,\cdots,\nabla h_k$ linearly independent.
By H\"ormander's equivalence of phase functions theorem
(\cite{Hormander}) we can microlocally rewrite $I_\psi$
in the form
\begin{equation}\label{EQ:new-Ipsi}
 I_\psi f(x)=\int_\Rn\p{\int_{\R^k} e^{i\sum_{j=1}^k \lambda_j
 h_j(x,y)} a(x,\overline{\lambda}) f(y) d\overline{\lambda}}
 dy,
\end{equation}
where $\overline{\lambda}=(\lambda_1,\ldots,\lambda_k)$ and
$a\in S^{(n-k)/2}_{1,0}(\Rn\times\R^{k})$ is (microlocally)
elliptic.
Now, let us take $f$ in the form $f=(I-\Delta)^{-s/2}\delta_{y_0}$
for some $y_0$ in the smooth part of the set
$$\Sigma^{y}=\b{y\in\Rn: (x,y)\in\Sigma \textrm{ for some }
 x}.$$
It follows that $f\in L^q_{loc}$ if and only if
$s>n(1-1/q)$. Now, let
$b(x,\overline{\lambda})\in
S^{-s+(n-k)/2}_{1,0}(\Rn\times\R^{k})$ be the amplitude of the
Fourier integral operator $I_\psi\circ (I-\Delta)^{-s}$.
Denoting $\overline{h}=(h_1,\ldots,h_k)$, we easily find that
$$I_\psi f(x)=(2\pi)^k \FT^{-1}_{\overline{\lambda}}b(x,
\overline{h}(x,y_0))\approx \abs{\textrm{dist}\,
(x,\Sigma_{y_0})}^{-k+s-(n-k)/2},$$
locally uniformly in $x$,
where $\Sigma_{y_0}$ is the set of all
$x\in\Rn$ such that $(x,y_0)\in\Sigma$.
Since $I_\psi f$ is smooth along
$\Sigma_{y_0}$, we find that
$I_\psi f\not\in L^q_{loc}(\Rn)$ if and only if
$s\leq k(1-1/q)+(n-k)/2$. Thus, if
$n(1-1/q)< k(1-1/q)+(n-k)/2$, operator $I_\psi$ is not locally
continuous on $L^q(\Rn)$. Since we assumed that $1\leq q<2$ and
that $I_\psi$ is locally continuous on $L^q(\Rn)$, it follows that
$k=0$, which means that $\psi$ must be linear.
\end{proof}

\begin{rem}\label{Rem:pert-lin}
Theorem \ref{TH:local-lq} has another interesting
relation with the $L^p$--properties of Fourier integral
operators. Suppose that the inverse
$\psi^{-1}$ of $\psi\in C^\infty(\Rn)$ can be written
in the form
$$
% \psi(x)=A x+\epsilon(x), \quad
 \psi^{-1}(x)=Ax+\delta(x),
$$
for some real-valued non-degenerate
matrix $A$ and
$\delta\in S^0(\Rn)$ with
$\|D\delta\|\ll 1$. Then using the
expression \eqref{EQ:alt-norms} for norms we have
$$
\|I_\psi f\|_{M^{p,q}_{loc}} =
 \left\| \|\Phi(D-k) I_\psi f\|_{L^p_{x,loc}}
 \right\|_{l^q_k}.
$$
Taking $f\in M^{p,q}_{comp}$, we have
\begin{align*}
  \|\Phi(D-k) I_\psi f\|_{L^p_{x,loc}} & =
  \| \Phi(D-k)\iint
     e^{i(x\cdot\psi^{-1}(\xi)-y\cdot\xi)} |\det D\psi^{-1}(\xi)|
     f(y) d\xi dy\|_{L^p_{x,loc}} \\
    & = \|\Phi(D-k) e^{ix\delta(D)} |\det D\psi^{-1}(D)|\,
    |\det A^{-1}|\, f\|_{L^p_{x,loc}}.
\end{align*}
Now, we observe the estimate
$$\|\Phi(D-k)e^{ix\delta(D)}g\|_{L^p_{loc}}=
\left\|\Phi(D-k)\sum_{j=0}^\infty\frac{(ix\delta(D))^j}{j!}
g(x)\right\|_{L^p_{loc}}\leq C\|\Phi(D-k)g\|_{L^p_{loc}},
$$
which holds since $x$ is bounded and $\delta$ is of order zero.
Combining these estimates together and using the boundedness of
$|\det D\psi^{-1}(D)|$ on $L^p$, we get that
$I_\psi$ is locally bounded on $M^{p,q}$ for $1<p<\infty$.
\end{rem}

Before proving Theorem \ref{homogen},
we first consider the simple case of it allowing a
harmonic analysis interpretation.

\par

\begin{rem}
Let us prove that
\begin{equation}\label{Hilbert-1}
\|f(|\cdot|)\|_{M^{p,q}}
\le C\|f\|_{M^{p,q}}
\quad \text{for all $f \in M^{p,q}(\R)$},
\end{equation}
where $1<p,q<\infty$.
Since
\begin{align*}
f(|x|)
&=f(|x|)\, \chi_{(-\infty,0)}(x)+f(|x|)\, \chi_{[0,\infty)}(x)
\\
&=f(-x)\, \chi_{(-\infty,0)}(x)+f(x)\, \chi_{[0,\infty)}(x),
\end{align*}
if $\MT_{\chi_{(-\infty,0)}}$ and $\MT_{\chi_{[0,\infty)}}$
are bounded on $M^{p,q}(\R^n)$, then
\begin{align*}
\|f(|\cdot|)\|_{M^{p,q}}
&\le \|\MT_{\chi_{(-\infty,0)}}(f(-\cdot))\|_{M^{p,q}}
+\|\MT_{\chi_{[0,\infty)}}f\|_{M^{p,q}}
\\
&\le \left(\|\MT_{\chi_{(-\infty,0)}}\|_{\calL(M^{p,q})}
+\|\MT_{\chi_{[0,\infty)}}\|_{\calL(M^{p,q})}\right) \|f\|_{M^{p,q}},
\end{align*}
that is, we obtain \eqref{Hilbert-1},
where $\MT_{\chi}$ is the operator of multiplication by $\chi$
(see Proposition \ref{prop:mult}).
Hence, it is enough to prove the boundedness of
$\MT_{\chi_{(-\infty,0)}}$ and $\MT_{\chi_{[0,\infty)}}$ on
$M^{p,q}(\R)$. By Proposition \ref{op-equiv},
if $\chi_{(-\infty,0)}(D)$ and $\chi_{[0,\infty)}(D)$
are bounded on $W^{q,p}(\R)$,
then $\MT_{\chi_{(-\infty,0)}}$ and $\MT_{\chi_{[0,\infty)}}$
are also bounded on $M^{p,q}(\R)$.
Let us prove the boundedness of
$\chi_{(-\infty,0)}(D)$ and $\chi_{[0,\infty)}(D)$ on $W^{p,q}(\R)$
for all $1<p,q<\infty$.
We recall that
\begin{equation}\label{EQ:alt-norms-2}
\|f\|_{W^{p,q}(\R)}
\asymp
\left\|\left( \sum_{k \in \Z}
|\Phi(D-k)f|^q\right)^{1/q}\right\|_{L^p(\R)},
\end{equation}
(see the proof of Theorem \ref{TH:mw}).
On the other hand,
it is known that,
for all $1<p,q<\infty$,
\begin{equation}\label{Hilbert-2}
\left\|\left(\sum_k |Hf_k|^q\right)^{1/q}\right\|_{L^p(\R)}
\le C_{p,q}\left\|\left(\sum_k |f_k|^q\right)^{1/q}\right\|_{L^p(\R)},
\end{equation}
where $H$ is the Hilbert transform,
\[
Hf(x)=\calF^{-1}[(-i\, \mathrm{sgn}\, \xi)\widehat{f}],
\qquad \mathrm{sgn}\, \xi=
\begin{cases}
1, &\xi>0
\\
0, &\xi=0
\\
-1, &\xi<0.
\end{cases}
\]
(see \cite[Theorem 8.1]{Duo}).
Since
\[
\chi_{(-\infty,0)}(\xi)=-(\mathrm{sgn}\, \xi-1)/2
\quad \text{and} \quad
\chi_{[0,\infty)}(\xi)=(\mathrm{sgn}\, \xi+1)/2
\]
for all $\xi \neq 0$, we have
\begin{equation}\label{Hilbert-3}
\chi_{(-\infty,0)}(D)=-(iH-I)/2
\quad \text{and} \quad
\chi_{[0,\infty)}(D)=(iH+I)/2
\end{equation}
where $If=f$.
Combining \eqref{EQ:alt-norms-2},
\eqref{Hilbert-2} and \eqref{Hilbert-3},
we see that  $\chi_{(-\infty,0)}(D)$ and $\chi_{[0,\infty)}(D)$
are bounded on $W^{p,q}(\R)$.
\end{rem}

\par

In the general case of Theorem \ref{homogen},
the proof is based on some investigations
of Gabor expansions of elements in $M^{p,q}$ and $W^{p,q}$.
More precisely, let $\{ x_j \} _{j\in I}$ and $\{ \xi _k \}_{k\in I}$
be lattices in $\rr n$, and consider functions or distributions of the form
$$
f(x)=\sum _{j,k\in I}c_{j,k}e^{i\scal x{\xi _k}}\chi (x-x_j),
$$
for some sequences $c=\{ c_{j,k}\} _{j,k\in I}$ and
$\chi \in M^{p_0}\setminus 0$, where $1\le p_0\le 2$.
We note that $f$ makes sense as an element in $M^{p_0}$
when $c$ belongs to $l^1_0$,
the set of all sequences $d=\{ d_{j,k}\} _{j,k\in I}$ such that
$d_{j,k}=0$, except for a finite numbers of $j$ and $k$.
We are especially concerned with finding conditions
on $p\in [1,\infty]$ and $q\in [1,\infty ]$ such that
$f$ still makes sense when $c$ belongs to $l_1^{p,q}$ or $l_2^{p,q}$.
Here $l_1^{p,q}$ consists of all sequences
$d=\{ d_{j,k}\} _{j,k\in I}$ such that
$$
\nm d{l_1^{p,q}}\equiv \Big ( \sum _{k\in I} \Big (
\sum _{j\in I} |c_{j,k}|^p\Big )^{q/p}\Big )^{1/q}<\infty
$$
(with obvious interpretation when $p=\infty$ or $q=\infty$), and $l_2^{p,q}$
consists of all sequences $d=\{ d_{j,k}\} _{j,k\in I}$ such that
%$$
%\nm d{l_2^{p,q}}\equiv \Big ( \sum _{j\in I}
%\Big ( \sum _{k\in I} |c_{j,k}|^p\Big )^{q/p}\Big )^{1/q}<\infty
%$$
$$
\nm d{l_2^{p,q}}\equiv \Big ( \sum _{j\in I}
\Big ( \sum _{k\in I} |c_{j,k}|^q\Big )^{p/q}\Big )^{1/p}<\infty
$$

\par

We have the following proposition.
%Here and in what follows we let $p'\in [1,\infty ]$
%be the conjugate exponent of $p\in [1,\infty ]$, i.{\,}e. $1/p+1/p'=1$.

\par

\begin{prop}\label{expansions1}
Assume that $p,p_0,q\in[1,\infty ]$ satisfy $1\le p_0\le \min (p,p',q,q')$,
and let  $\{ x_j \} _{j\in I}$ and $\{ \xi _k \}_{k\in I}$
be lattices in $\rr n$. Then the map
$$
(\, \{ c_{j,k}\} _{j,k\in I} ,\chi )\mapsto
\sum _{j,k\in I}c_{j,k}e^{i\scal x{\xi _k}}\chi (x-x_j)
$$
from $l^1_0\times M^{1}$ to $M^{1}$ extends uniquely to a continuous map
from $l_1^{p,q}\times M^{p_0}$ to $M^{p,q}$,
and from $l_2^{p,q}\times M^{p_0}$ to $W^{p,q}$.
Furthermore, for some constant $C$ it holds
\begin{align}
\Big \Vert  \sum _{j,k\in I}c_{j,k}e^{i\scal x{\xi _k}}\chi (x-x_j)
\Big \Vert _{M^{p,q}}
&\le C\nm {\{ c_{j,k}\} _{j,k\in I}}{l_1^{p,q}}\nm \chi {M^{p_0}}
\label{normsum1}
\intertext{and}
\Big \Vert  \sum _{j,k\in I}c_{j,k}e^{i\scal x{\xi _k}}\chi (x-x_j)
\Big \Vert _{W^{p,q}}
&\le C\nm {\{ c_{j,k}\} _{j,k\in I}}{l_2^{p,q}}\nm \chi {M^{p_0}}.
\label{normsum2}
\end{align}
\end{prop}

\par

\begin{proof}
We only prove the mapping properties for $l_1^{p,q}$.
The other case follows by similar arguments and is left for the reader.

\par

We may assume that $p_0= \min (p,p',q,q')$. First we observe that the
result holds in the case $p_0=1$ or $p_0=2$,
in view of the general Gabor theory on modulation spaces
(cf. Chapters 6 and 12 in \cite{Grochenig}).
The result now follows for general $p_0$, by multi-linear interpolation
between these cases. The proof is complete.
\end{proof}

\par

Next we consider multiplication properties of $M^{p,q}$ spaces and $W^{p,q}$
with certain types of step functions.
For this reason it is convenient to make the following definition.

\par

\begin{defn}\label{funktionsrum}
Let
\begin{enumerate}
\item $\Sigma _0(\rr n)$ be the set of all functions $\uppsi $ such that
$$
\uppsi =\sum _{j\in \Lambda} c_j\chi _{x_j+Q}
$$
for some cube $Q\subseteq \rr n$, lattice $\{ x_j \} _{j\in \Lambda}
\subseteq \rr n$ and sequence $\{ c_j \} _{j\in \Lambda} \in
l^\infty$;

\vrum

\item $\Sigma (\rr n)$ be the set of all functions $\uppsi $ such that
$$
\uppsi =\sum _{j\in \Lambda} \fy _j\chi _{x_j+Q}
$$
for some cube $Q\subseteq \rr n$, lattice $\{ x_j \} _{j\in \Lambda}
\subseteq \rr n$ and sequence $\{ \fy _j \} _{j\in \Lambda} \subseteq
C^\infty (\rr n)$ such that $\{ \partial ^\alpha \fy _j \} _{j\in
\Lambda}$ is a bounded sequence in $L^\infty $ for every multi-index
$\alpha$.
\end{enumerate}
\end{defn}

\par

We have now the following result.

\par

\begin{prop}\label{thm1}
Assume that $\psi \in \Sigma (\mathbb R^n)$ and $1<p,q<\infty$. Then the map
$M_\psi $ from $\mathscr S(\mathbb R^n)$
to $\mathscr S'(\mathbb R^n)$ extends uniquely to a continuous map on
$M^{p,q}(\mathbb R^n)$ and on $W^{p,q}(\mathbb R^n)$.
\end{prop}
%{\bf modification for $p,q=\infty?$

\par

\begin{proof}
It is no restriction to assume that $\uppsi$ is as in Definition 
\ref{funktionsrum}
with $\{ x_j \} _{j\in J}=\mathbb Z^n$ and
$$
Q=\sets {x\in \rr n}{0\le x_1,\dots ,x_n\le 1}.
$$
We only prove the assertion for $M^{p,q}$.
The other case follows by similar arguments and is left for the reader.

\par

First we assume that $\psi \in \Sigma _0$,
and we let $\chi _j(t)$ for $j=0,1,2$ and $t\in \mathbb R$ be defined by the 
formulas
$$
\chi _0(t)=\max (1-|t|,0),\quad \chi _1(t)=\chi _0(t)\chi _{(-1,0)}(t),\quad
\text{and}\quad  \chi _2(t)=\chi _0(t)\chi _{(0,1)}(t),
$$
where $\chi _{(a,b)}$ is the characteristic function of the interval $(a,b)$.
By straight-forward computations it follows that $\chi _0\in M^1(\mathbb R)$,
and that $\chi _{(-1,0)},\chi _{(0,1)}\in M^{1,q}$ for every $q>1$. Hence 
Proposition \ref{prop:mult} gives
$$
\chi _1,\chi _2 \in M^1\cdot M^{1,p_0}\subseteq M^{\infty ,1}\cdot M^{1,p_0}
\subseteq M^{1,p_0}\subseteq M^{p_0},
$$
when $p_0= \min (p,p',q,q')>1$, and if
$$
\kappa _0\equiv \chi _0\otimes \cdots \otimes \chi _0,\quad
\kappa _l\equiv \chi _{l_1}\otimes \cdots \otimes \chi _{l_n},\quad
l=(l_1,\dots ,l_n)\in \{ 1,2\} ^n,
$$
with $n$ factors in the tensor products, then $\kappa _0\in M^1(\rr n)$, $\kappa 
_l\in M^{p_0}(\rr n)$
when $l\in \{ 1,2\} ^n$, and
\begin{equation}\label{superpos1}
\kappa _0=\sum _{l}\kappa _l.
\end{equation}

\par

Next assume that $f\in M^{p,q}$ is arbitrary,
and let $\{ x_j\} _{j\in J}=\mathbb Z^n$. Then
$$
f(x)=\sum _{j,k\in J}c_{j,k}e^{i\scal x{\xi _k}}\kappa _0(x-x_j),
$$
for some lattice $\{ \xi _k \} _{k\in J}$ and sequence
$\{ c_{j,k}\} \in l_1^{p,q}$.
This follows from the general Gabor theory for modulation spaces
(cf. Chapter 12 in \cite{Grochenig}). Furthermore, by
Proposition \ref{expansions1} it follows that
$$
f_l(x)\equiv \sum _{j,k\in J}c_{j,k}e^{i\scal x{\xi _k}}\kappa _l(x-x_j)
$$
makes sense as an element in $M^{p,q}$ for every $l\in \{ 1,2\} ^n$,
and, hence
$$
f=\sum _lf_l
$$
in view of \eqref{superpos1}. It therefore suffices to prove that
$\uppsi \cdot f_l\in M^{p,q}$ for every $l$.

\par

First assume that $\{ c_{j,k}\} \in l^1_0$. From the assumptions if
follows that
$$
\uppsi (x) =\sum _{j\in J}d_j\chi _Q(x-x_j),
$$
where $\chi _Q$ is the characteristic function of the unit cube
$[0,1]^n$, and $\{ d_j\}  \in l^\infty$.

Then $\uppsi \cdot f_l\in M^{p,q}$
is well-defined, and by the definitions we have
$$
\uppsi (x)\cdot f_l(x)=
\sum _{j,k\in J}\widetilde c_{j,k}e^{i\scal x{\xi _k}}\kappa
_l(x-x_j),
$$
where
$$
\widetilde c_{j,k} = c_{j,k}d_{j}.
$$

\par

Since $\{ d_j\} \in l^\infty$, it follows that
$$
\nm {\{ \widetilde c_{j,k} \}}{l_1^{p,q}}\le  \nm {\{ c_{j,k}
\}}{l_1^{p,q}}\nm {\{ d_j\} }{l^\infty}<\infty .
$$
Hence $\uppsi f_l\in M^{p,q}$ in view of Proposition \ref{expansions1},
and the result follows in this case.

\par

Next assume that $\uppsi \in \Sigma$ is arbitrary, $\fy _j$ as in
Definition \ref{funktionsrum}{\,}(2), and let $C>0$. Then we may split
up $\{ x_j\} _{j\in \Lambda}$ into sublattices $\{ x_j\} _{j\in \Lambda
_1}$,\dots ,$\{ x_j\} _{j\in \Lambda _N}$ such
that if $j_1,j_2\in \Lambda _m$ and $j_1\neq j_2$ for some $1\le m\le
N$, then the distance $d_{j_1,j_2}$ between $x_{j_1}+Q$ and
$x_{j_2}+Q$ is larger than $C$. Now set
$$
\uppsi _m = \sum _{j\in \Lambda _m}\fy _j\chi _{x_j+Q}.
$$
Since $\uppsi =\sum _m\uppsi _m$, the result follows
if we prove that the map $f\mapsto \uppsi _m\cdot f$ extends
uniquely to a continuous map on $M^{p,q}$ and on $W^{p,q}$
for every $1\le m\le N$.

\par

From the fact that $d_{j_1,j_2}\ge C$ it follows that there is a
non-negative function $\phi \in C_0^\infty (\rr n)$ such that $\phi =1$
on $Q$ and $\supp \phi _{j_1}\cap \supp \phi _{j_2}=\emptyset$ when
$j_1,j_2\in \Lambda _m$ for some $m$ and $j_1\neq j_2$. Here $\phi
_j=\phi (\cdot -x_j)$. This gives $\uppsi _m=b\, c$ where
$$
b = \sum _{j\in \Lambda _m}\chi _{x_j+Q}\quad \text{and}\quad c = \sum
_{j\in \Lambda _m} \phi _j\fy _j
$$
In particular, $c$ is a smooth function on $\rr n$ and bounded
together with all its derivatives, and $b\in \Sigma _0(\rr n)$,
which imply that
$$
c\in M^{\infty ,1}\subseteq W^{\infty ,1}.
$$
By the first part of the proof, and Proposition \ref{prop:mult}, it
follows that $\mathscr M_b$ and
$\mathscr M_c$ are bounded on $M^{p,q}$ and on $W^{p,q}$.

\par

Hence
$$
\nm {\uppsi _m\cdot f}{M^{p,q}} = \nm {\mathscr M_b\mathscr M_c
f}{M^{p,q}}\le C_1\nm f{M^{p,q}},\quad f\in \mathscr S,
$$
for some constant $C_1$, and similarly
when the $M^{p,q}$ norms are replaced by $W^{p,q}$ norms.
This proves that $\mathscr M_\psi$ extends to
continuous mappings on $M^{p,q}$ and on $W^{p,q}$. It also follows
that these extensions are unique since $\mathscr S$ is dense in
$M^{p,q}$ and $W^{p,q}$. The proof is complete.
\end{proof}

\par

\begin{proof}[Proof of Theorem \ref{homogen}]
Assume that $f\in \mathscr S$. For any $\theta \in \{ 0,1\} ^n$, set
$$
g_\theta (x_1,\dots ,x_n) = f (S(x)+T((-1)^{\theta _1}x_1,\dots
,(-1)^{\theta _n}x_n))\chi ((-1)^{\theta _1}x_1,\dots ,(-1)^{\theta
_n}x_n),
$$
where $\chi$ is the characteristic function of the set
$$
\sets {x\in \rr n}{x_j>0}.
$$
Since compositions by affine mappings are continuous operations on
modulation spaces, Proposition \ref{thm1} shows that the map $f\mapsto
g_\theta$ from  $\mathscr S$ to $\mathscr S'$ is uniquely extendable
to a continuous map on $M^{p,q}$ and on $W^{p,q}$. The assertions is
now a consequence of the fact that
$$
f(S(x)+T(|x_1|,\dots ,|x_n|)) =
\sum _{\theta \in \{ 0,1\} ^n} g_\theta,\qquad \text{in
$\mathscr S'({\rr n})$},
$$
when $f\in \mathscr S({\rr n})$. The proof is complete.
\end{proof}

\section{Wiener type spaces on the torus, and properties of periodic
distributions}

In this section we will indicate counterparts of our results
in the case of spaces on the torus as well as make several
related observations. Some of these remarks concern Wiener type
properties for periodic distributions. In fact, there is a
one-to-one corresponding between periodic distributions (periodic
continuous functions) and distributions (continuous functions),
respectively, on the torus.

\par

We fix the notation
$\Tn=(\R/2\pi\Z)^n$ as well as the Fourier transform and its
inverse given by
$$
  (\FT f)(\xi) = \widehat{f}(\xi)
  = (2\pi)^{-n}\int_{\Tn} e^{-ix\cdot\xi} f(x) dx,
\quad
  f(x) = \sum_{\xi\in\Z^n} e^{ix\cdot\xi}\widehat{f}(\xi).
$$
For the analysis of pseudo-differential operators on the
torus using
the Fourier series and for the justification of operators
below we refer to \cite{RT07}.
We note also that canonical transforms on
the torus can be viewed as a special case of
Fourier series operators considered in \cite{RT07}.

A straightforward modification of the definitions of the modulation
and Wiener amalgam spaces from \eqref{EQ:alt-norms} is
\begin{align*}\label{EQ:alt-norms-torus}
\|f\|_{M^{p,q}(\Tn)}\asymp\left\| \|\Phi(D-k)f\|_{L^p_x(\Tn)}
\right\|_{l^q_k(\Z^n)}, \\
\|f\|_{W^{p,q}(\Tn)}\asymp\left\| \|\Phi(D-k)f\|_{l^q_k(\Z^n)}
\right\|_{L^p_x(\Tn)},
\end{align*}
for some $\Phi$ with compact support in $\Z^n$ (in the
discreet topology).
Then we can easily observe the equalities (which can be regarded
as a counterpart of Theorem \ref{TH:mw} in the local setting)
\begin{equation}\label{EQ:eq-spaces-torus}
M^{p,q}(\Tn)=W^{p,q}(\Tn)=\FT l^q(\Tn) \textrm{ for all }
1\leq p,q\leq\infty.
\end{equation}
In particular, in the case of $q=1$ this equality can be 
viewed as a characterisation of 
absolutely convergent Fourier series.

Further, as a counterpart of Theorem \ref{TH:B-H}
the Beurling--Helson property automatically holds
on $M^{p,q}(\Tn)$ and $W^{p,q}(\Tn)$ because of the
original Beurling and Helson theorem \cite{BH} ($q=1$) as well
as extensions for other $q$ (\cite{LO94, Se76}).

As a counterpart of Theorem \ref{TH:local-lq} we observe that if
$\psi:\Z^n\to\Z^n$ is a bijection then the canonical transform
$I_\psi$ is bounded on $M^{p,q}(\Tn)$
for all $1\leq p,q\leq \infty$.
Indeed, let $f\in M^{p,q}(\Tn)$.
Then $\widehat{f}\in l^q(\Z^n)$ and hence
$$
\n{I_\psi f}_{M^{p,q}(\Tn)}=
\n{\psi^*\widehat{f}}_{l^q(\Z^n)}=
\p{\sum_{\xi\in\Z^n} \abs{\widehat{f}(\psi(\xi))}^q}^{1/q}=
\n{\widehat{f}}_{l^q(\Z^n)}\asymp
\n{f}_{M^{p,q}(\Tn)}.
$$
Next, we recall the result of \cite{BGOR} that
for $0\leq \alpha\leq 2$
operators $e^{i|D|^\alpha}$ are bounded on modulation
spaces $\calM^{p,q}(\Rn)$, $1\leq p,q\leq\infty$
(for the definition see Remark \ref{modified-def}).
In particular, this covers wave and Schr\"odinger propagators.

To give a counterpart of Remark \ref{rem:propagators}
on the torus,
using \eqref{EQ:eq-spaces-torus} we easily
conclude that propagators $e^{i|D|^\alpha}$ are bounded on
$M^{p,q}(\Tn)$ and $W^{p,q}(\Tn)$ for all $1\leq p,q\leq\infty$
and all $\alpha\in\R$. Moreover, they are isometries
on these spaces if we induce their norms from the space
$\FT l^q(\Tn)$.

We also note that results of this section can be extended
to general compact Lie groups $G$ if we use a natural extension
of the global definition of modulation and Wiener amalgam spaces
using the duality between $G$ and the space of its 
continuous irreducible
unitary representations $\widehat{G}$ (as in \cite{RT-book}).

\medspace

We finish this section with
the following proposition for periodic distributions,
parallel to \eqref{EQ:eq-spaces-torus}. Here we refer to Appendix A for the
definition of the weighted spaces $M^{p,q}_{(\omega )}$ and
$W^{p,q}_{(\omega )}$.

\par

\begin{prop}
Assume that $p,q\in [1,\infty ]$, $f\in \mathscr S'(\rr n)$ is
periodic, and that $\omega \in \mathscr P(\rr {2n})$ is such that
$\omega (x,\xi )=\omega _0(\xi )$, for some $\omega _0\in \mathscr
P(\rr n)$. Then the following conditions are equivalent:
\begin{enumerate}
\item $f\in M^{\infty ,q}_{(\omega )}(\rr n)$;

\vrum

\item $\chi \cdot f\in M^{p,q}_{(\omega )}(\rr n)$, for each $\chi \in
\mathscr S(\rr n)$;

\vrum

\item $f\in W^{\infty ,q}_{(\omega )}(\rr n)$;

\vrum

\item $\chi \cdot f\in W^{p,q}_{(\omega )}(\rr n)$, for each $\chi \in
\mathscr S(\rr n)$;

\vrum

\item $\chi \cdot f\in \mathscr FL^{q}_{(\omega _0)}(\rr n)$, for each
$\chi \in \mathscr S(\rr n)$.
\end{enumerate}
\end{prop}

\par

\begin{proof}
Since $\mathscr S(\rr n)$ is contained in each space of the form
$M^{p,q}_{(v)}(\rr n)$ and $W^{p,q}_{(v)}(\rr n)$, when $v\in \mathscr
P(\rr {2n})$, it follows from Proposition A.1 in Appendix A that (1)
implies (2) and (3) implies (4).

\par

Next assume that (2) is fulfilled, and consider $F(x,\xi )=V_\fy
f(x,\xi )$, where $\fy \in C_0^\infty \setminus 0$. Then
$F(x,\xi )$ is a smooth function and period in the $x$-variable,
with the same period $t\in \rr n$ as $f$. Let $Q\subseteq \rr n$ be a
cube with side length $t$, and let $\chi \in C_0^\infty (\rr n)$ be
equal to $1$ in the set $Q+\supp \fy$. Then we have
\begin{multline*}
\nm f{M^{\infty ,q}_{(\omega )}} =\Big ( \int \underset {x \in \R^n}
{\operatorname{ess\, sup}}|V_\fy f(x,\xi )|^q\, d\xi \Big )^{1/q}
\\[1ex]
=\Big ( \int \underset {x \in Q}
{\operatorname{ess\, sup}}|V_\fy f(x,\xi )|^q\, d\xi \Big )^{1/q}
=\Big ( \int \underset {x \in Q}
{\operatorname{ess\, sup}}|(V_\fy (\chi \cdot f))(x,\xi )|^q\, d\xi
\Big )^{1/q}
\\[1ex]
=\nm {\chi f}{M^{\infty ,q}_{(\omega )}}\le C\nm {\chi
f}{M^{p,q}_{(\omega )}}<\infty .
\end{multline*}
This proves that (2) is equivalent to (1).

\par

In the same way it follows that (3) is equivalent to (4). In
particular, if (2) or (4) are fulfilled for a particular $p$, then
they are fulfilled for any $p\in [1,\infty ]$. Hence,
(A.1) in Appendix A gives that (2), (4) and (5) are equivalent. This
proves the result.
\end{proof}

\par

%%%%%%%%%%%%%%%%%%%%%%%%%%%%%%%%%%
\section*{Appendix A, Some remarks on weighted Wiener type spaces}
%%%%%%%%%%%%%%%%%%%%%%%%%%%%%%%%%%

In this appendix we make some reviews of general (or weighted)
modulation spaces and weighted versions of $W^{p,q}$, and some
multiplication and convolution proerties of such spaces. We start to
consider appropriate conditions on the involved weight functions.

\par

Assume that $0<\omega ,v\in L^\infty _{loc}(\rr n)$. Then $\omega$ is
called $v$-moderate, if $\omega (x+y)\le C\omega (x)v(y)$, for some
constant $C$ which is independent of $x,y\in \rr n$. If in addition,
$v$ can be chosen as a polynomial, then $\omega$ is called polynomial
moderated. We let $\mathscr P(\rr n)$ be the set of polynomial
moderated functions on $\rr n$. For any $\omega \in \mathscr P(\rr
n)$, it follows that
$$
P(x)^{-1}\le \omega (x)\le P(x),
$$
for some polynomial $P$ on $\rr n$.

\par

Next assume that $\varphi \in \calS(\R^n)\setminus 0$ and $\omega
\in
\mathscr P(\rr {2n})$ are fixed. Then the modulation space
$M^{p,q}_{(\omega )}(\rr n)$ consists of all $f\in \mathscr S'(\rr n)$
such that $V_{\fy}f(x,\xi )\omega (x,\xi )\in L^{p,q}_1(\rr {2n})$,
and with equipp by the norm
$\nm f{M^{p,q}_{(\omega )}}\equiv \nm {V_\fy f\, \omega}{L^{p,q}_1}$
is finite. We also let the Wiener amalgam related space
$W^{p,q}_{(\omega )}(\rr n)$ be the set of $f\in \mathscr S'(\rr n)$
such that $\nm f{M^{p,q}_{(\omega )}}\equiv \nm {V_\fy f\,
\omega}{L^{p,q}_2}$ is finite.

\par

If $\omega \in \mathscr P(\rr n)$, then we let $L^p_{(\omega
)}(\rr n)$ be the set of all $f\in L^1_{loc}(\rr n)$ such that $\nm
{f\, \omega}{L^p}<\infty$, and we let $\mathscr FL^q_{(\omega )}(\rr
n)$ be the set of all $f\in \mathscr S'(\rr n)$ such that $\nm
{\widehat f\, \omega}{L^q}<\infty$. By Proposition 1.7 in \cite{To04}
and Theorem 3.2 in \cite{To04B} it follows that the embeddings
\begin{equation}\tag*{(A.1)}
\begin{alignedat}{2}
M^{p,q_1}_{(\omega )}(\rr n) &\subseteq L^p_{(\omega _0)}(\rr
n)\subseteq M^{p,q_2}_{(\omega )}(\rr n),&\qquad \omega (x,\xi )
&=\omega _0(x )\in \mathscr P(\rr n),
\\[1ex]
M^{p_1,q}_{(\omega )}(\rr n) &\subseteq \mathscr FL^q_{(\omega
_0)}(\rr n)\subseteq M^{p_2,q}_{(\omega )}(\rr n),&\qquad \omega
(x,\xi ) &= \omega _0(\xi )\in \mathscr P(\rr n),
\\[1ex]
W^{p,q_1}_{(\omega )}(\rr n) &\subseteq L^p_{(\omega _0)}(\rr
n)\subseteq W^{p,q_2}_{(\omega )}(\rr n),&\qquad \omega (x,\xi ) &=
\omega _0(x)\in \mathscr P(\rr n),
\\[1ex]
W^{p_1,q}_{(\omega )}(\rr n) &\subseteq \mathscr FL^q_{(\omega
_0)}(\rr n)\subseteq W^{p_2,q}_{(\omega )}(\rr n),&\qquad \omega
(x,\xi ) &= \omega _0(\xi )\in \mathscr P(\rr n),
\end{alignedat}
\end{equation}
hold for each $p,p_j,q,q_j\in [1,\infty ]$ for $j=1,2$ such that
$$
p_1\le \min (q,q'),\quad p_2\ge \max (q,q'),\quad q_1\le \min
(p,p'),\quad q_2\ge \max (p,p').
$$

\par

Almost all properties for non-weighted modulation spaces and Wiener
amalgam spaces can be generalised to spaces of the form
$M^{p,q}_{(\omega )}$ and $W^{p,q}_{(\omega )}$. For example these
spaces are Banach spaces, and independent of the choice of window
function $\fy \in \mathscr S(\rr n)\setminus 0$, where different
choices of $\fy$ give rise to equivalent norms. Furthermore, if
$p_1\le p_2$, $q_1\le q_2$ and $\omega _2\le C\omega _1$ for some
constant $C$, then $M^{p_1,q_1}_{(\omega _1)}\subseteq
M^{p_2,q_2}_{(\omega _2)}$. We also have that
\begin{equation}\tag*{(\ref{(2.1)})$'$}
 \|f\|_{W^{p,q}_{(\omega )}}\asymp\|\widehat{f}\|_{M^{q,p}_{(\omega
_0)}},\qquad \omega _0(\xi ,-x)=\omega (x,\xi ),
\end{equation}
and we note that
\begin{alignat*}{2}
M^{p,q}_{(\omega )}(\rr n) &\subseteq W^{p,q}_{(\omega )}(\rr n)&\quad
\text{when}\quad q\le p
\intertext{and}
W^{p,q}_{(\omega )}(\rr n) &\subseteq M^{p,q}_{(\omega )}(\rr n)&\quad
\text{when}\quad p\le q.
\end{alignat*}

\medspace

Next we discuss multiplication and convolution properties for
modulation spaces. Assume that $\omega _0,\dots ,\omega _N\in \mathscr
P(\rr {2n})$, $p_0,\dots ,p_N\in [1,\infty ]$ and $q_0,\dots ,q_N\in
[1,\infty ]$ satisfy
\begin{equation}\tag*{(A.2)}  %%\label{weightyoung1}
\begin{aligned}
\omega _0(x,\xi _1+\cdots +\xi _N)&\le C\omega _1(x,\xi _1)\cdots
\omega _N(x,\xi _N),
\\[1ex]
\frac 1{p_1}+\cdots +\frac 1{p_N} &= \frac 1{p_0}\quad \text{and}
\quad \frac
1{q_1}+\cdots +\frac 1{q_N} = N-1+\frac 1{q_0} ,
\end{aligned}
\end{equation}
for some constant $C$ which is independent of $x,\xi _1,\dots ,\xi
_N\in \rr n$. Then
\begin{align}
\nm {f_1\cdots f_N}{M^{p_0,q_0}_{(\omega _0)}} &\le C^N\nm
{f_1}{M^{p_1,q_1}_{(\omega _1 )}}\cdots \nm {f_N}{M^{p_N,q_N}_{(\omega
_N)}}\tag*{(A.3)}
\intertext{and}
\nm {f_1\cdots f_N}{W^{p_0,q_0}_{(\omega _0)}} &\le C^N\nm
{f_1}{W^{p_1,q_1}_{(\omega _1 )}}\cdots \nm {f_N}{W^{p_N,q_N}_{(\omega
_N)}},\tag*{(A.4)}
\end{align}
for some constant $C$ which is independent of $N$ and $f_1,\dots
f_N\in \mathscr S(\rr n)$. Here the first inequality is a consequence
of \cite[Theorem 3]{Fe3} and its proof. The second inequality is an
immediate consequence of \cite[Theorem 5.5]{To04B}. By
Hahn-Banach's theorem it follows that the map
$$
(f_1,\dots ,f_N)\mapsto f_1\cdots f_N
$$
from $\mathscr S\times \cdots \times \mathscr S$ to
$\mathscr S$ extends to a continuous multiplication from
$M^{p_1,q_1}_{(\omega _1)}\times \cdots \times M^{p_N,q_N}_{(\omega
_N)}$ to $M^{p_0,q_0}_{(\omega _0)}$, and from
$W^{p_1,q_1}_{(\omega _1)}\times \cdots \times W^{p_N,q_N}_{(\omega
_N)}$ to $W^{p_0,q_0}_{(\omega _0)}$.

\par

If instead $\omega _0,\dots ,\omega _N\in \mathscr
P(\rr {2n})$, $p_0,\dots ,p_N\in [1,\infty ]$ and $q_0,\dots
,q_N\in [1,\infty ]$ satisfy
\begin{equation}\tag*{(A.5)}
\begin{aligned}
\omega _0(x_1+\cdots +x_N,\xi )&\le C\omega _1(x_1,\xi )\cdots \omega
_N(x _N,\xi ),
\\[1ex]
\frac 1{p_1}+\cdots +\frac 1{p_N} &= N-1+\frac 1{p_0}\quad
\text{and}\quad \frac 1{q_1}+\cdots +\frac 1{q_N} = \frac 1{q_0},
\end{aligned}
\end{equation}
for some constant $C$ which is independent of $x_1,\dots ,x_N,\xi \in
\rr n$, then we have
\begin{align}
\nm {f_1*\cdots *f_N}{M^{p_0,q_0}_{(\omega _0)}} &\le C^N\nm
{f_1}{M^{p_1,q_1}_{(\omega _1 )}}\cdots \nm {f_N}{M^{p_N,q_N}_{(\omega
_N)}},\tag*{(A.6)}
\\[1ex]
\nm {f_1*\cdots *f_N}{W^{p_0,q_0}_{(\omega _0)}} &\le C^N\nm
{f_1}{W^{p_1,q_1}_{(\omega _1 )}}\cdots \nm {f_N}{W^{p_N,q_N}_{(\omega
_N)}},\tag*{(A.7)}
\end{align}
for some constant $C$ which is independent of $N$ and $f_1,\dots
f_N\in \mathscr S(\rr n)$. By Hahn-Banach's theorem it follows that
the convolution map
$$
(f_1,\dots ,f_N)\mapsto f_1*\cdots *f_N
$$
from $\mathscr S\times \cdots \times \mathscr S$ to
$\mathscr S$ extends to a continuous multiplication from
$M^{p_1,q_1}_{(\omega _1)}\times \cdots \times M^{p_N,q_N}_{(\omega
_N)}$ to $M^{p,q}_{(\omega )}$, and from
$W^{p_1,q_1}_{(\omega _1)}\times \cdots \times W^{p_N,q_N}_{(\omega
_N)}$ to $W^{p,q}_{(\omega )}$.

\par

A problem here concerns the \emph{uniqueness} for the extensions of
multiplications and convolutions, since it easily appears that
there may be
situations where more than one of those $p_j$ or $q_j$ are allowed to
be equal to $\infty$. Consequently, $\mathscr S$ might fail to be dense
in more than one of the involved modulation or Wiener amalgam related
spaces. In these situations, we define multiplications and
convolutions between elements in modulation spaces in the same way as
in \cite {To04,To04B}, using the formulae
\begin{equation}\tag*{(A.8)}
\begin{aligned}
&({f_1\cdots  f_N},{g}) = \iint \big ( V_{\fy _1}f_1(x ,\cdo )*\cdots
*V_{\fy _N}f_N(x ,\cdo  )\big )(\xi )\overline {V_{\fy _0}g(x,\xi
)}{\,}dxd\xi ,
\\[1ex]
&\text{where}\ \, \fy _0,\dots ,\fy _N\in \mathscr S(\rr n) \ \,
\text{satisfy}\ \,  \int \fy _1(x)\cdots \fy _N(x)\overline {\fy
_0(x)}\, dx =(2\pi )^{-Nn}.
\end{aligned}
\end{equation}
and
\begin{equation}\tag*{(A.9)}
\begin{aligned}
&({f_1*\cdots * f_N},{\fy}) = \iint \big ( V_{\fy _1}f_1(\cdo ,\xi
)*\cdots *V_{\fy _N}f_N(\cdo ,\xi )\big )(x)\overline {V_{\fy
_0}g(x,\xi )}{\,}dxd\xi ,
\\[1ex]
&\text{where}\ \, \fy _0,\dots ,\fy _N\in \mathscr S(\rr n) \ \,
\text{satisfy}\ \,  \int (\fy _1*\cdots *\fy _N)(x)\overline {\fy
_0(x)}\, dx =(2\pi )^{-n},
\end{aligned}
\end{equation}
when $f_1,\dots ,f_N,g\in \mathscr S(\rr n)$ (cf. (2.3) in \cite{To04}
and (5.4) in \cite{To04B}). Theorem 5.5 in \cite{To04B} and its proof
then shows that the following propositions are true:

\par

\renewcommand{\rubrik}{Proposition A.1}

\begin{tom}
Assume that $p_j,q_j\in [1,\infty ]$ and $\omega _j\in \mathscr P(\rr
{2n})$ for $j=0,\dots ,N$ satisfy {\rm{(A.2)}} for some constant $C$,
independent of $x, \xi _1,\dots ,\xi _N \in \rr n$. Then the following
is true:
\begin{enumerate}
\item $(f_1,\dots ,f_N)\mapsto f_1\cdots f_N$ is a continuous, symmetric
and associative map from $M^{p_1,q_1}_{(\omega _1)}(\rr n)\times
\cdots \times M^{p_N,q_N}_{(\omega _N)}(\rr n)$ to $M^{p_0,q_0}
_{(\omega _0)} (\rr n)$, which is independent of the choice of $\fy
_0,\dots ,\fy _N$ in {\rm{(A.8)}}. Furthermore, {\rm{(A.3)}} holds for
some constant $C$ which is independent of $f_j\in M^{p_j,q_j}_{(\omega
_j)}(\rr n)$ for $j=1,\dots ,N$;

\vrum

\item $(f_1,\dots ,f_N)\mapsto f_1\cdots f_N$ is a continuous, symmetric
and associative map from $W^{p_1,q_1}_{(\omega _1)}(\rr n)\times
\cdots \times W^{p_N,q_N}_{(\omega _N)}(\rr n)$ to $W^{p_0,q_0}
_{(\omega _0)} (\rr n)$, which is independent of the choice of $\fy
_0,\dots ,\fy _N$ in {\rm{(A.8)}}. Furthermore, {\rm{(A.4)}} holds for
some constant $C$ which is independent of $f_j\in M^{p_j,q_j}_{(\omega
_j)}(\rr n)$ for $j=1,\dots ,N$.
\end{enumerate}
\end{tom}

\par

\renewcommand{\rubrik}{Proposition A.2}

\begin{tom}
Assume that $p_j,q_j\in [1,\infty ]$ and $\omega _j\in \mathscr P(\rr
{2n})$ for $j=0,\dots ,N$ satisfy {\rm{(A.5)}} for some constant $C$,
independent of $x_1,\dots ,x_N, \xi \in \rr n$. Then the following is
true:
\begin{enumerate}
\item $(f_1,\dots ,f_N)\mapsto f_1*\cdots *f_N$ is a continuous,
symmetric and associative map from $M^{p_1,q_1}_{(\omega _1)}(\rr
n)\times \cdots \times M^{p_N,q_N}_{(\omega _N)}(\rr n)$ to $M^{p_0,q_0}
_{(\omega _0)} (\rr n)$, which is independent of the choice of $\fy
_0,\dots ,\fy _N$ in {\rm{(A.9)}}. Furthermore, {\rm{(A.6)}} holds for
some constant $C$ which is independent of $f_j\in M^{p_j,q_j}_{(\omega
_j)}(\rr n)$ for $j=1,\dots ,N$;

\vrum

\item $(f_1,\dots ,f_N)\mapsto f_1*\cdots *f_N$ is a continuous, symmetric
and associative map from $W^{p_1,q_1}_{(\omega _1)}(\rr n)\times
\cdots \times W^{p_N,q_N}_{(\omega _N)}(\rr n)$ to $W^{p_0,q_0}
_{(\omega _0)} (\rr n)$, which is independent of the choice of $\fy
_0,\dots ,\fy _N$ in {\rm{(A.9)}}. Furthermore, {\rm{(A.7)}} holds for
some constant $C$ which is independent of $f_j\in M^{p_j,q_j}_{(\omega
_j)}(\rr n)$ for $j=1,\dots ,N$.
\end{enumerate}
\end{tom}

%%%==================================================================
%%%==================================================================

\end{document}